\documentclass[11pt]{article}
\usepackage{latexsym}
\newcommand{\bbz}{\mbox{\boldmath $Z$}}

\newcommand{\bbc}{\mbox{\boldmath $C$}}
\newcommand{\bbq}{\mbox{\boldmath $Q$}}
\newcommand{\bbf}{{F}}
\newcommand{\bbr}{\mbox{\boldmath $R$}}
\newcommand{\qed}{\hfill$\Box$}
\newcommand{\too}{\longrightarrow}
\newcommand{\then}{\Longrightarrow}
\newcommand{\G}{{X}}
\newcommand{\C}{{\cal K}}
\newcommand{\T}{{\rho}}
\newcommand{\sn}{{\it span}}
\newcommand{\tr}{{\it Tr}}

\newcommand{\D}{{\it Des}}

\newcommand{\al}{\alpha}

\newcommand{\tw}{{\tilde{w}}}
\newcommand{\ts}{{\tilde{s}}}
\newcommand{\bs}{{\dot{s}}}
\newcommand{\hgt}{{ht}}

\newcommand{\ua}{{\dot{a}}}
\newcommand{\ub}{{\dot{b}}}
\newcommand{\ud}{{\dot{d}}}
\newcommand{\da}{{\ddot{a}}}
\newcommand{\db}{{\ddot{b}}}
\newcommand{\dd}{{\ddot{d}}}
\newcommand{\oa}{{\vec a}}
\newcommand{\ob}{{\vec b}}

\newtheorem{thm}{Theorem}[section]
\newtheorem{pro}[thm]{Proposition}
\newtheorem{lem}[thm]{Lemma}

\newtheorem{con}[thm]{Conjecture}
\newtheorem{cor}[thm]{Corollary}
\newtheorem{fac}[thm]{Fact}
\newtheorem{prb}[thm]{Problem}

\newtheorem{obs}[thm]{Observation}
\newtheorem{rem}[thm]{Remark}
\newtheorem{df}[thm]{Definition}

\begin{document}

\title{A Unified Construction of Coxeter Group Representations}
\bibliographystyle{acm}
\author{Ron M. Adin%
\thanks{Department of Mathematics, Bar-Ilan University,
Ramat-Gan 52900, Israel. Email: {\tt radin@math.biu.ac.il}}\ $^\S$
\and Francesco Brenti%
\thanks{Dipartimento di Matematica,
Universit\'{a} di Roma ``Tor Vergata'',
Via della Ricerca Scientifica,
00133 Roma, Italy. Email: {\tt brenti@mat.uniroma2.it}}\ $^\S$
\and Yuval Roichman%
\thanks{Department of Mathematics, Bar-Ilan University,
Ramat-Gan 52900, Israel. Email: {\tt yuvalr@math.biu.ac.il}}
\thanks{Research of all authors was supported in part by
the Israel Science Foundation, founded by the Israel Academy of Sciences and Humanities;
by EC's IHRP Programme, within the Research Training Network ``Algebraic Combinatorics
in Europe'', grant HPRN-CT-2001-00272;
and by internal research grants from Bar-Ilan University.}}
\date{February 15, 2006 (revised)}


\maketitle

\begin {abstract}
An elementary approach to the construction of Coxeter group
representations is presented.
\end{abstract}




\section{Introduction}\label{s.intro}

\subsection{Outline}\label{s.out}

Let $(W,S)$ be a Coxeter system, and let $\C$ be a subset of $W$.
Let $\bbf$ be a suitable field of characteristic zero (e.g., the
field $\bbc(q)$ in the case of the Iwahori-Hecke algebra), and let
$\T$ be a representation of (the Iwahori-Hecke algebra of) $W$ on
the vector space $V_{\C}:=\sn_{\bbf}\{C_w\,|\,w\in\C\}$, with
basis vectors indexed by elements of $\C$.
We study the sets $\C$ and representations $\T$ which satisfy the
following axiom:

\begin{itemize}
\item[$(A)$] {\it For any generator $s\in S$ and any element $w\in
\C$ there exist scalars $a_s(w),b_s(w)\in\bbf$ such that
$$
\T_s(C_w) = a_{s}(w)C_w + b_{s}(w)C_{ws}.
$$
If $w\in\C$ but $ws\not\in\C$ we assume $b_s(w)=0$.%
}
\end{itemize}

A pair $(\T,\C)$ satisfying Axiom $(A)$ is called an {\em abstract
Young (AY) pair;} $\T$ is an {\em AY representation,} and $\C$ is
an {\em AY cell.} If $\C\ne \emptyset$ and has no proper subset
$\emptyset \subset \C' \subset \C$ such that $V_{\C'}$ is
$\T$-invariant, then $(\T,\C)$ is called a {\em minimal AY pair.}
(This is much weaker than assuming $\T$ to be irreducible.)

\bigskip

Surprisingly, Axiom $(A)$ leads to very concrete matrices, whose
entries are essentially inverse linear. Analysis of the
construction involves a convexity theorem of Tits~\cite{T} and the
generalized descent classes introduced by Bj\"orner and
Wachs~\cite{BW}.

\medskip

This paper has been motivated by an attempt to extend the classical 
representation theory of the symmetric group~\cite{Ja,JK,Md,Sa,St2} to 
arbitrary Coxeter groups. Our method of generalization is through the 
explicit representation matrices, following ideas from~\cite{BR,OV,V}. 
Having completed an early version of this paper~\cite{U1}, we have been
informed of the recent papers~\cite{Ram03, Cher}, which give an
essentially equivalent construction.
This paper presents a more elementary approach. The implications
of Axiom (A) are studied in detail. In particular, the Coxeter relations
are studied and it is shown that, under mild conditions,
the diagonal coefficients in Axiom (A) determine all the character
values (solving a problem presented by Kazhdan).
It is also shown that, under the same mild conditions, the
coefficients in Axiom (A) are determined by a linear functional, the
one that appears as a basic ingredient in~\cite{Ram03}. The paper is
also focused on combinatorial aspects; in particular, the
combinatorial structure of induced AY representations is studied
in Section~\ref{s.induce}.

\bigskip

The paper is organized as follows. The rest of this section lists
our main results. Section~\ref{s.prelim} contains preliminaries.
The main concept of this paper, abstract Young pairs, is
introduced in Section~\ref{s.AY_pairs}. In
Sections~\ref{s.AY_pairs}--\ref{s.CRC}
we discuss general Coxeter groups. The basic axiom on
representation matrices is analyzed in Sections~\ref{s.AY_pairs},
\ref{s.reflect} and Appendix 1 (Section~\ref{s.ax1}). Minimal AY
cells and their connection to convexity are studied in
Section~\ref{s.convexity}. Implications of Coxeter relations with
$m(s,t)\le 3$ (within general Coxeter groups) are discussed in
Section~\ref{s.CRC}. In Section~\ref{s.sl} it is shown that, for
simply laced Coxeter groups, minimal AY representations are
determined by a linear functional. The Hecke algebra analogues are
studied in Section~\ref{s.sl_q}. AY cells associated to induced
representations are studied in Section~\ref{s.induce}. Examples of
AY representations are given in Section~\ref{examples}. Finally,
the independence of character values on non-diagonal coefficients
is proved in Appendix 2 (Section~\ref{s.ax2}).


\subsection{Main Results: Cells}\label{s.main_cells}

Recall the definition of AY cells and representations from
the previous subsection.

\begin{prb}\label{p.rep-cell1}{\rm (Kazhdan~\cite{Ka})}
Given a subset $\C\subseteq W$, how many nonequivalent
abstract Young representations may be defined on $V_\C$?
\end{prb}
In particular,
\begin{prb}\label{p.rep-cell2}
Which subsets of $W$ are (minimal) AY cells?
\end{prb}

In this paper we give a partial solution to Problem~\ref{p.rep-cell2};
see Section~\ref{s.mincells}.
In particular, it will be shown that a minimal AY cell must be convex.

\begin{pro}\label{main.convex}{\rm (see Corollary~\ref{t.MAY_cell}(i))}
Every minimal AY cell is convex
(in the right Cayley graph $\G(W,S)$ or, equivalently,
under right weak Bruhat order).
\end{pro}

\begin{obs}\label{t.main_unitcell}{\rm (see Observation~\ref{t.unitcell})}
Every nonempty AY cell is a left translate of an AY cell containing
the identity element of $W$.
\end{obs}

AY cells in the symmetric group will be characterized
in~\cite{U-II}.



A combinatorial rule for induction of AY representations is given
in Section~\ref{s.induce}. This rule is analogous to the one for
Kazhdan-Lusztig representations, proved by Barbasch and Vogan for
Weyl groups~\cite[Prop. 3.15]{BV2} and by Geck for Coxeter groups
\cite{Geck}.

\begin{thm}\label{t.main.i4} (see Theorem~\ref{t.i4})
Let $(W,S)$ be a finite Coxeter system, $P=\langle J\rangle$
$(J\subseteq S)$ a parabolic subgroup,
and $W^J$ the set of all representatives of minimal length of the
right cosets of $P$ in $W$. Let $(\psi,\cal D)$ be a minimal AY
pair for $P$. Then
\begin{itemize}
\item[1.] ${\cal D} W^J$ is
a minimal AY cell for $W$. \item[2.] The induced representation
$\psi\uparrow ^W_P$ is isomorphic to an AY representation on
$V_{{\cal D} W^J}$.
\end{itemize}
\end{thm}

In the case of the symmetric group, the cells associated to induced AY
representations consist of shuffles. This result is analogous
to~\cite[K Rule V]{GR}; see also~\cite{Sc}.

\subsection{Main Results: Representations}\label{s.main_reps}

In Section~\ref{s.reflect} and Appendix 1 it is shown that,
under mild conditions, Axiom $(A)$ is equivalent to the following
more specific version.
Here $T$ is the set of all reflections in $W$.

\begin{itemize}
\item[$(B)$]
{\it
For any reflection $t\in T$
there exist scalars $\ua_t, \ub_t, \da_t, \db_t \in\bbf$
such that, for all $s\in S$ and $w\in \C$:
$$
\T_s(C_w) =\cases{%
\ua_{wsw^{-1}}C_w+\ub_{wsw^{-1}}C_{ws}, &if $\ell(w)<\ell(ws)$;\cr
\da_{wsw^{-1}}C_w+\db_{wsw^{-1}}C_{ws}, &if $\ell(w)>\ell(ws)$.\cr}
$$
If $w\in\C$ and $ws\not\in\C$ we assume that
$\ub_{wsw^{-1}}=0$ (if $\ell(w)<\ell(ws)$) or
$\db_{wsw^{-1}}=0$ (if $\ell(w)>\ell(ws)$).%
}
\end{itemize}

\begin{thm}\label{main.a1}{\rm (see Theorem~\ref{t.a1})}
Let $(\T,\C)$ be a minimal AY pair for the Iwahori-Hecke algebra of $(W,S)$.
If 
$a_s(w)=a_{s'}(w')\then b_s(w)=b_{s'}(w')$
$(\forall s,s'\in S, w,w'\in\C)$,
then $\rho$ satisfies Axiom $(B)$.
\end{thm}

This theorem shows that the coefficients $a_s(w)$ and $b_s(w)$ in Axiom $(A)$
depend only on the reflection $wsw^{-1}\in T$ and on the relation between
$w$ and $ws$ in right weak Bruhat order.

The assumption regarding the coefficients $b_s(w)$ in Theorem~\ref{main.a1}
is merely a normalization condition;
for more details see Subsection~\ref{s.normal}.

Thus, in order to determine an AY representation, it suffices to determine
the coefficients $\ua_t$ for all reflections $t$
($\da_t$ is determined by $\ua_t$, see Lemma~\ref{t.l1}(a))
and to choose a normalization for the $\ub_t$.

\begin{prb}\label{p.char}{\rm (Kazhdan~\cite{Ka})}
Do the coefficients $\ua_t$ determine the character values?
\end{prb}

An affirmative answer to this problem is given in Appendix 2:

\begin{thm}\label{t.b_indep.main} {\rm (see Theorem~\ref{t.b_indep})}
Let $(\T,\C)$ be an AY pair satisfying Axiom $(B)$.
Then all the corresponding character values are polynomials in
$a$-coefficients only (no $b$-coefficients).
\end{thm}

By Observation~\ref{t.main_unitcell}, every abstract Young representation is
isomorphic to one on an AY cell containing the identity element.
Therefore, in the rest of this subsection, we assume $id\in\C$.

It turns out that for simply laced Coxeter groups the coefficients $\ua_t$
are given by a linear functional.


\begin{thm}\label{t.main11}{\rm (See Theorem~\ref{t.f2rep})}
Let $(W,S)$ be an irreducible simply laced Coxeter system,
and let $\C$ be a convex subset of $W$ containing the identity element.
Let $\langle\ ,\ \rangle$ be an arbitrary positive definite
bilinear form on the root space $V$.
If $f\in V$ is $\C$-generic then
$$
\ua_{wsw^{-1}} := {1\over \langle f,\al_{wsw^{-1}}\rangle}
\qquad(\forall w\in \C, s\in S),
$$
together with $\da_{wsw^{-1}}$, $\ub_{wsw^{-1}}$ and $\db_{wsw^{-1}}$
satisfying conditions (a) and (b) of Lemma~\ref{t.l1} (with $q = 1$),
defines a representation $\T$ such that $(\T,\C)$ is a minimal AY pair.
\end{thm}

\noindent
For the definition of $\C$-genericity see Definition~\ref{d.cg} below.

Replacing $\langle f,\al_t\rangle$
by its $q$-analogue $[\langle f,\al_t\rangle]_q$
gives representations of the Iwahori-Hecke algebra ${\cal H}_q(W)$.
See Theorem~\ref{t.f2rep.q} below.

\smallskip

The following theorem is complementary.

\begin{thm}\label{t.main12}{\rm (See Theorem \ref{t.rep2f})}
Let $(W,S)$ be an irreducible simply laced Coxeter system
and let $\C$ be a subset of $W$ containing the identity element.
If $(\T,\C)$ is a minimal AY pair satisfying Axiom $(B)$
and $\ua_{wsw^{-1}}\ne 0$ ($\forall w\in \C, s\in S$)
then
there exists a $\C$-generic $f\in V$ such that
$$
\ua_{wsw^{-1}}={1\over \langle f,\al_{wsw^{-1}}\rangle}\qquad
(\forall\ w\in \C, s\in S).
$$
\end{thm}

For an Iwahori-Hecke algebra analogue see Theorem~\ref{t.rep2f.q} below.

\bigskip


\bigskip

%

A combinatorial bijection between elements of minimal AY cells and
standard Young tableaux of the same shape will be given
in~\cite{U-II}.
This is used to prove that

\begin{thm}\label{main.irr} {\rm \cite{U-II}}
Every irreducible representation of the symmetric group $S_n$
may be realized as a minimal abstract Young representation.
\end{thm}

This result is extended to the hyperoctahedral groups $B_n$
in~\cite{U-II}. More examples of AY representations are given in
Section~\ref{examples}.



\section{Preliminaries}\label{s.prelim}

A {\em Coxeter system} is a pair $(W,S)$ consisting of a group $W$ and
a set $S$ of generators for $W$, subject only to relations of the form
$$
(st)^{m(s,t)}=1,
$$
where $m(s,s)=1$ and $m(s,t)=m(t,s)\ge 2$ for $s\ne t$ in $S$.
In case no relation occurs for a pair $(s,t)$, we make the convention
that $m(s,t)=\infty$. $W$ is called a {\em Coxeter group}.
If $m(s,t)\le 3$ for all $s\ne t$ then $(W,S)$ is called {\em simply laced}.
Throughout this paper we assume that the set $S$ of generators is finite;
the group $W$ itself may be infinite.

The {\em Iwahori-Hecke algebra} (sometimes called the {\em generic algebra})
${\cal H}_q(W)$ of a Coxeter system $(W,S)$
is the algebra generated by $\{\hat{s}\,|\,s\in S\}$
with the involution relation on generators replaced by
\begin{equation}\label{e.deformed_inv}
(\hat{s}-1)(\hat{s}+q_s)=0 \qquad (\forall s\in S),
\end{equation}
in addition to the usual braid relations: for any two distinct generators
$s,t\in S$
$$
\hat{s} \hat{t} \hat{s} \cdots=
\hat{t} \hat{s} \hat{t} \cdots \qquad
(m(s,t) \hbox{ factors on each side}).
$$
Here $q_s$ is an indeterminate depending only on the conjugacy class (in $W$)
of the simple reflection $s$.
For convenience, we shall use the notation $q_t:=q_s$ for any reflection $t\in T$;
here $s\in S$ is any generator (simple reflection) conjugate to $t$ in $W$.
For example, if $W$ is an irreducible simply laced
Coxeter group then all reflections are conjugate, and there is only one
indeterminate $q$. If $W$ is an irreducible finite Weyl group
then there are at most two values of $q$, corresponding to
the classification of the roots into long and short.

Note that our conventions are slightly non-standard;
for the standard version, replace in (\ref{e.deformed_inv}) above
each generator $\hat{s}$ by $-\hat{s}$.

Let $(W,S)$ be a Coxeter system, and let $\G=\G(W,S)$ be the corresponding
Cayley graph (with generators acting on the right):
its vertices are the elements of $W$, and $x, y\in W$ are connected by an edge
if and only if $x^{-1}y\in S$. $\G$ is a connected undirected graph.
A subset $\C \subseteq W$ is called {\em convex} if,
for any $x, y\in \C$, all the geodesics (paths of shortest length)
connecting $x$ to $y$ in $\G$ have all their vertices in $\C$.

Let $P$ be a poset, and let $\G$ be its undirected Hasse diagram.
Thus $\G$ is an undirected graph which has $P$ as a vertex set,
with an edge $\{x,y\}$ whenever $x$ either covers or is covered by $y$.
A subset $\C$ of $P$ is called {\em convex} if, for any $x,y\in\C$,
all geodesics (shortest paths) connecting $x$ to $y$ in $\G$ have
all their vertices in $\C$.

In this paper, $P$ will be a Coxeter group $W$ with
the right weak Bruhat order, namely the transitive closure of the relation
$$
w<ws \iff w\in W, s\in S \hbox{\ and\ } \ell(w)<\ell(ws).
$$
Clearly, $\C\subseteq W$ is convex in the right Cayley graph $\G(W,S)$
if and only if it is convex in the right weak Bruhat poset $P$.

\bigskip

Let $V$ be the root space of a Coxeter system $(W,S)$.
Thus $V$ is a vector space over $\bbr$ with a basis $\{\al_s\,|\,s\in S\}$
indexed by the group generators.
A symmetric bilinear form $B$ is defined on $V$ by
$$
B(\al_s,\al_t) := -\cos{\pi\over m(s,t)}\qquad(\forall s,t\in S)
$$
(interpreted to be $-1$ in case $m(s,t)=\infty$).
For each generator $s\in S$ define a linear map $\sigma_s:V \to V$ by
$$
\sigma_s(v) := v - 2B(v,\al_s)\al_s\qquad(\forall v\in V).
$$
This yields a faithful $B$-preserving action $\sigma$ of $W$ on $V$
(see, e,g.,~\cite{Hum}),
and defines the corresponding {\em root system}
$$
\Phi := \{\sigma_w(\al_s)\,|\,w\in W, s\in S\} \subseteq V.
$$
Let $T$ be the set of all conjugates of generators;
its elements are called {\em reflections} of $W$,
and elements of $S$ are called {\em simple reflections}.
There is a bijection between reflections in $T$ and {\em positive roots}
in $\Phi$, given by
$t=wsw^{-1} \longleftrightarrow \al_t=\sigma_w(\al_s)$
(provided that $\ell(w)<\ell(ws)$).

Throughout the paper $\langle \ ,\ \rangle$ is an arbitrary positive definite
bilinear form on $V$.



The {\em height} of a positive root $\al_t=\sum\limits_{s\in S}c_s
\al_s$ is the positive integer
$$
\hgt(\al_t) := \sum\limits_{s\in S} c_s.
$$

%
%
%

\section{Abstract Young Pairs}\label{s.AY_pairs}

In this section we set the general axiomatic framework of the paper.

Throughout the paper we shall fix a base field $\bbf$ of characteristic zero.
When assuming the symmetric-orthogonal normalization
(see Subsection~\ref{s.normal} below),
we may take $\bbf = \bbr$ or $\bbf = \bbc$.
When assuming the stochastic normalization, we may also take $\bbf = \bbq$ .
For the corresponding Hecke algebra, we should take a suitable field of
rational functions in $q$ (or the various $q_s$).

Let $(W,S)$ be a Coxeter system, and let $\C$ be a subset of $W$.
Let $\bbf$ be a suitable field of characteristic zero as above,
and let $\T$ be a representation of (the Iwahori-Hecke algebra of) $W$
on the vector space
$$
V_{\C}:=\sn_{\bbf}\{C_w\,|\,w\in\C\},
$$
with basis vectors indexed by the elements of $\C$.
We shall study the sets $\C$ and representations $\T$ which satisfy
the following axiom:

\begin{itemize}
\item[$(A)$]
{\it
For any generator $s\in S$ and any element $w\in \C$
there exist scalars $a_s(w),b_s(w)\in\bbf$ such that
$$
\T_s(C_w) = a_{s}(w)C_w + b_{s}(w)C_{ws}.
$$
If $w\in\C$ but $ws\not\in\C$ we assume $b_s(w)=0$.%
}
\end{itemize}

\begin{df}
A pair $(\T,\C)$ satisfying Axiom $(A)$ is called an
{\em abstract Young (AY) pair;} $\T$ is an {\em AY representation,}
and $\C$ is an {\em AY cell.}
If $\C$ is nonempty and has no proper subset
$\emptyset \subset \C' \subset \C$ such that $V_{\C'}$ is $\T$-invariant,
then $(\T,\C)$ is called a {\em minimal AY pair,}
and $\T$ (respectively, $\C$) is
a {\em minimal AY cell} (respectively, {\em representation}).
\end{df}

Axiom $(A)$ may be stated more abstractly. Indeed, consider the following
statement:

\begin{itemize}
\item[$(A')$]
{\it
For any parabolic subgroup $H=\langle J\rangle$ of $W$ (where $J\subseteq S$)
and any $w\in W$, the ``coset subspace''
$V_{wH \cap \C}$
is invariant under $\T(h)$, for all $h\in H$.%
}
\end{itemize}

\begin{obs}\label{t.aprime}
Axiom $(A)$ is the special case $|J|=1$ (i.e., $H$ minimal parabolic)
of Axiom $(A')$. Actually, it is easy to verify that Axioms $(A)$ and $(A')$
are equivalent.
\end{obs}

\begin{obs}\label{t.unitcell}
Left multiplication by a group element is an automorphism of the right
Cayley graph $\G(W,S)$ (see Section~\ref{s.prelim}).
Therefore, if $\C\subseteq W$ carries a representation $\T$ satisfying
Axiom $(A)$ then, for any $v\in W$, the left translate $v^{-1}\C$ carries
an equivalent representation (with exactly the same representation matrices).
\end{obs}

\section{Convexity}\label{s.convexity}

\subsection{Minimal Cells and Convexity}\label{s.mincells}

In this subsection we give a partial solution (Corollary~\ref{t.MAY_cell}(i))
to Problem~\ref{p.rep-cell2}.

\begin{df}
Let $(\T,\C)$ be an AY pair.
An arc $v \to vs$ in the Cayley graph $\G(W,S)$ is {\em $\T$-feasible} if
$v,vs\in \C$ and
$b_s(v)\ne 0$.
A path $v \to vs_1 \to \ldots \to vs_1 \cdots s_k$ is {\em $\T$-feasible} if
all its arcs are $\T$-feasible, i.e., if
$$
b_{s_i}(vs_1 \cdots s_{i-1})\ne 0 \qquad(1\le i\le k).
$$
If there is a $\T$-feasible path from $v$ to $w$, denote
$$
v \Rightarrow_\T w.
$$
\end{df}

%

Recall that a {\em geodesic} from $v$ to $w$ in $\G(W,S)$ is a path
from $v$ to $w$ of shortest possible length (i.e., length $\ell(v^{-1}w)$).

\begin{lem}\label{t.feas_geo}
Let $(\T,\C)$ be an AY pair, and let $v,w\in\C$.
If there exists a $\T$-feasible geodesic (in $\C$) from $v$ to $w$,
then any geodesic in $X(W,S)$ from $v$ to $w$ is $\T$-feasible.
\end{lem}
{\bf Proof.}
Let $v \to vs_1 \to \ldots \to vs_1 \cdots s_k = w$ be a $\T$-feasible geodesic.
Then $s_1 \cdots s_k$ is a reduced word for $v^{-1}w$, and
$b_{s_i}(vs_1 \cdots s_{i-1})\ne 0$ ($1\le i\le k$).
By assumption,
$vs_1 \cdots s_i \in\C$ ($0\le i\le k$).
Consider the vector $\T_{s_k} \cdots \T_{s_1}(C_v) \in V_{\C}$.
The support of this vector contains only basis vectors of the form
$C_{v\ts_1 \cdots \ts_k}$, where each $\ts_i$ is either $s_i$ or the identity
element of $W$. Since $s_1 \cdots s_k$ is reduced,
$C_{v\ts_1 \cdots \ts_k} = C_w$ only if
$\ts_i = s_i$ for all $i$, and the
coefficient of $C_w$ in $\T_{s_k} \cdots \T_{s_1}(C_v)$ is therefore
$$
b_{s_k}(vs_1 \cdots s_{k-1}) \cdots b_{s_1}(v) \ne 0.
$$
On the other hand,
$\T_{s_k} \cdots \T_{s_1} = \T_{s_k \cdots s_1} = \T_{(v^{-1}w)^{-1}}$
is independent of the choice of reduced word for $v^{-1}w$.
Thus, if
$v \to vt_1 \to \ldots \to vt_1 \cdots t_k = w$
is any geodesic from $v$ to $w$ then
$t_1 \cdots t_k$ is a reduced word for $v^{-1}w$, and
$C_w$ is also in the support of $\T_{t_k} \cdots \T_{t_1}(C_v)$.
Similar reasoning now shows that the coefficient of $C_w$ in
$\T_{t_k} \cdots \T_{t_1}(C_v)$ is
$b_{t_k}(vt_1 \cdots t_{k-1}) \cdots b_{t_1}(v)$, and therefore
$$
b_{t_i}(vt_1 \cdots t_{i-1})\ne 0 \qquad(1\le i\le k).
$$
Thus the geodesic $v \to vt_1 \to \ldots \to vt_1 \cdots t_k = w$
is $\T$-feasible.

\qed

\begin{lem}\label{t.feas_path_geo}
Let $(\T,\C)$ be an AY pair, and let $v,w\in\C$.
If there exists a $\T$-feasible path (of any length) from $v$ to $w$,
then there exists a $\T$-feasible geodesic from $v$ to $w$.
\end{lem}
{\bf Proof.}
Let $v \to vs_1 \to \ldots \to vs_1 \cdots s_k = w$ be the shortest
$\T$-feasible path from $v$ to $w$.
If $k = \ell(v^{-1}w)$ then the path is a geodesic, and we are done.

Assume now that $k > \ell(v^{-1}w)$.
Let $j:= \max\{i\,|\,\ell(s_1 \cdots s_i) = i\}$. Then $1\le j\le k-1$.
Denote $u:= vs_1 \cdots s_j$, so that $\ell(v^{-1}u) = j$ whereas
$\ell(v^{-1}us_{j+1}) < j+1$, hence
$\ell(v^{-1}us_{j+1}) = j-1$.
Let $t_1 \cdots t_{j-1}$ be a reduced word for $v^{-1}us_{j+1}$, so that
$t_1 \cdots t_{j-1}s_{j+1}$ is a reduced word for $v^{-1}u$.
The path $v \to vs_1 \to \ldots \to vs_1 \cdots s_k$ is $\T$-feasible, and
therefore its initial segment $v \to vs_1 \to \ldots \to vs_1 \cdots s_j = u$
is a $\T$-feasible geodesic.
$v \to vt_1 \to \ldots \to vt_1 \cdots t_{j-1}
\to vt_1 \cdots t_{j-1}s_{j+1} = u$
is another geodesic connecting the same pair of vertices,
so by Lemma~\ref{t.feas_geo} it is $\T$-feasible.
Consider the path
$$
v \to vt_1 \to \ldots \to vt_1 \cdots t_{j-1} = vs_1 \cdots s_{j+1} \to
vs_1 \cdots s_{j+2} \to \ldots \to vs_1 \cdots s_k = w.
$$
It is a concatenation of two $\T$-feasible paths, and is therefore
$\T$-feasible. It has length $k-2$, contradicting our choice of $k$ as minimal.
Thus indeed $k=\ell(v^{-1}w)$, and the proof is complete.

\qed

Recall (from Section~\ref{s.prelim}) the definition of convexity in
the right Cayley graph, or equivalently in the right weak Bruhat poset.

\begin{cor}\label{t.MAY_cell}
Let $(\T,\C)$ be a minimal AY pair. Then:
\begin{itemize}
\item[{\bf (i)}]
$\C$ is convex;
\item[{\bf (ii)}]
$s\in S, w, ws\in\C \then b_s(w)\ne 0.$
\end{itemize}
\end{cor}
{\bf Proof.}
For $v\in\C$, let $\C(v)$ be the set of all $w\in\C$ for which there exists
a $\T$-feasible path from $v$ to $w$. Clearly, $\C(v)\ne\emptyset$
(since $v\in \C(v)$) and $V_{\C(v)}$ is $\T$-invariant. By minimality of $\C$,
$\C(v) = \C$. Thus, for any $v,w\in\C$, there exists a $\T$-feasible path
(and, by Lemma~\ref{t.feas_path_geo}, a $\T$-feasible geodesic)
from $v$ to $w$. By Lemma~\ref{t.feas_geo}, every geodesic from $v$ to $w$
is $\T$-feasible, and thus
its vertices
belong to $\C$. This shows that $\C$ is convex; the other claim follows by
taking a path (geodesic) of length one.

\qed

\begin{df}
An AY pair $(\T,\C)$ is called {\em strongly connected} if,
for any $v,w\in\C$,  $v \Rightarrow_\T w$
(i.e., there exists a $\T$-feasible path from $v$ to $w$).
\end{df}

\begin{cor}\label{t.conn_min}
An AY pair is strongly connected if and only if it is minimal.
\end{cor}

Thus, minimal AY cells are equivalence classes with respect to
the relation $\Rightarrow_\T$.
In this sense, the notion  of minimal AY cells
is similar to the notion of Kazhdan-Lusztig cells \cite{KL}.

\subsection{Convex Sets and Generalized Descent Classes}\label{s.convex_gdc}

Having seen the relevance of convex sets to our discussion,
we recall in this subsection, following Tits, some properties of convex sets.
In particular, we give an explicit combinatorial interpretation
to the equivalence relation in a minimal AY cell;
see Corollary~\ref{t.equiv_gdc} below.


Let $T$ be the set of all reflections in $W$, and let $A\subseteq T$
be any subset. The (left) {\em $A$-descent set} of an element $w\in W$
is defined by
$$
\D_A(w):=\{t\in A\,|\,\ell(tw)<\ell(w)\}.
$$
For $D\subseteq A\subseteq T$, the corresponding
{\em generalized descent class} is
$$
W_A^D := \{w\in W\,|\,\D_A(w)=D\}.
$$

\medskip

These sets were studied by Tits \cite[Ch. 2]{T} and Bj\"orner-Wachs
\cite{BW, BW2}.
Our notation follows \cite{BW}.

\smallskip

\noindent{\bf Examples.}
\begin{itemize}
\item[{\bf (i)}]
If $A=\emptyset$ then there is only one generalized descent class:
$W^\emptyset_\emptyset=W$.
\item[{\bf (ii)}]
If $A=S$ then the $A$-descent sets (and generalized descent classes) are
the (left) standard ones.
\item[{\bf (iii)}]
If $A=T$
then $\D_T(w)$ is the set of (left) associated reflections of $w$
(see, e.g., ~\cite[Ch. 3]{BB}). In particular, if
$W$ is the symmetric group then $\D_T(w)$ is the set of
(left) {\em inversions} of $w$.
\item[{\bf (iv)}]
$W_A^D$ may be empty; e.g., for $A=T$ and $D=S$ when $W$ is not commutative
(since $S\subseteq \D_T(w)$ if and only if $w$ is the longest element $w_0$
in the (necessarily finite) Coxeter group $W$, and then $\D_T(w_0)=T\ne S$).
\end{itemize}

The following theorem is a consequence of \cite[Theorem 2.19 ]{T}.
See \cite[Theorem 5.3]{BW} and \cite[Proposition 6.2]{BW2}.

\begin{thm}\label{t.tits} {\rm (Tits)}
A subset $\C$ of $W$ is convex under right weak Bruhat order
if and only if $\C=W_A^D$ for suitable subsets
$D\subseteq A\subseteq T$.
\end{thm}

We shall also use the following result.

\begin{pro}\label{t.des_diff} {\rm \cite[Prop. 3.1.3]{BB}}
If $w\in W$, $s\in S$ and $\ell(w)<\ell(ws)$ then
$$
\D_T(ws)=
\D_T(w)
\cup \{wsw^{-1}\}\qquad\mbox{\rm(disjoint union)}.
$$
\end{pro}

\begin{df}\label{d.A_equiv}
Let $A\subseteq T$.
Define a symmetric relation $\leftrightarrow_A$ on $W$ by
$$
w \leftrightarrow_A ws \ \ \iff \ \ wsw^{-1}\not\in A,
$$
and let $\sim_A$ be the transitive closure of $\leftrightarrow_A$.
The {\em $A$-cell} $\C_A(w)$ of an element $w\in W$ is the equivalence class
of $w$ under $\sim_A$.
\end{df}

\begin{thm}\label{t.thm3}
Let $u,v\in W$. Then
$$
u\sim_A v \ \ \iff \ \ \D_A(u)=\D_A(v).
$$
\end{thm}

\noindent{\bf Proof.}

Assume $u\sim_A v$. We may assume that $u \leftrightarrow_A v$.
Let $s \in S$ be such that $v=us$.
Then by Proposition~\ref{t.des_diff},
$$
\D_T(v)\setminus \{usu^{-1}\}=
\D_T(u)\setminus \{usu^{-1}\}.
$$
Since $u \leftrightarrow_A us$, $usu^{-1}\not\in A$ and therefore
$\D_A(v)=\D_A(u)$.

Conversely, assume that $\D_A(u)=\D_A(v)$.
Then $u,v\in W_A^D$ where $D:=\D_A(u)$. By Tits' Theorem, $W_A^D$ is convex.
Let $u,us_1,\dots, us_1\cdots s_r=v$ be a path of minimal length
connecting $u$ and $v$ in the Hasse diagram of the right weak order of $W$.
Then $us_1\cdots s_i\in W_A^D$ for all $0\le i\le r$, so
$\D_A(us_1\cdots s_i)=\D_{A}(us_{1} \cdots s_{i-1})$ for all $1\le i\le r$.
This implies that $t:=us_1\cdots s_{i-1}s_is_{i-1}\cdots s_1 u^{-1}\not\in A$,
and hence that
$us_1\cdots s_{i-1} \leftrightarrow_A us_1\cdots s_i$ for all $1\le i\le r$.

\qed

\begin{cor}\label{t.equiv_gdc}
For any $A\subseteq T$, the equivalence classes for the relation $\sim_A$
are exactly the nonempty generalized descent classes of the form $W_A^D$
($D\subseteq A$).
\end{cor}

\begin{cor}\label{t.convex_gdc}
For a subset $\C\subseteq W$, the following are equivalent:
\begin{itemize}
\item[{\bf (i)}]
$\C$ is convex.
\item[{\bf (ii)}]
$\C=W_A^D$ for suitable $D\subseteq A\subseteq T$.
\item[{\bf (iii)}]
$\C$ is empty, or $\C=\C_A(w)$
for some $w\in W$ and $A\subseteq T$.
\end{itemize}
(The subsets $A\subseteq T$ in {\bf(ii)} and {\bf(iii)} are the same.)
\end{cor}

\bigskip

By Corollary~\ref{t.MAY_cell}(i), every minimal AY cell is convex.
It should be noted that the converse does not hold: convex sets
are not necessarily minimal AY cells. Furthermore, for a given
$A\subseteq T$, $W_A^D$ may be a minimal AY cell for certain
$D\subseteq A$ but not for others. For example, let $W=S_5$,
$S=\{(i,i+1)\,|\,1\le i\le 4\}$ and
$A=\{(1,2),(2,3),(4,5),(1,4),(2,5)\}$. Then $\C_A(45123)$ is of
order $2$ and therefore not a minimal AY cell.
On the other hand, $\C_A(12345)$ is a minimal AY cell.
See~\cite{U-II}.

\section{Coefficients and Reflections}\label{s.reflect}

By Corollaries~\ref{t.MAY_cell}(i) and~\ref{t.convex_gdc},
every minimal AY cell is determined by certain sets of reflections.
Motivated by this observation, we will study the role of reflections
in determining the matrix entries for the corresponding AY representations.
This will lead to a substantial reformulation of Axiom $(A)$.

\subsection{Axiom $(B)$}\label{s.axiom_B}

Recall the definition of the set $T$ of all reflections in $W$.

Consider the following axiom for a pair $(\T,\C)$, where $\C$ is a subset of $W$
and $\T$ is a representation of (the Iwahori-Hecke algebra of) $W$ on the
vector space $V_{\C}:=\sn_{\bbf}\{C_w\,|\,w\in\C\}$:

\begin{itemize}
\item[$(B)$]
{\it
For any reflection $t\in T$
there exist scalars $\ua_t, \ub_t, \da_t, \db_t \in\bbf$
such that, for all $s\in S$ and $w\in \C$:
$$
\T_s(C_w) =\cases{%
\ua_{wsw^{-1}}C_w+\ub_{wsw^{-1}}C_{ws}, &if $\ell(w)<\ell(ws)$;\cr
\da_{wsw^{-1}}C_w+\db_{wsw^{-1}}C_{ws}, &if $\ell(w)>\ell(ws)$.\cr}
$$
If $w\in\C$ and $ws\not\in\C$ we assume that
$\ub_{wsw^{-1}}=0$ (if $\ell(w)<\ell(ws)$) or
$\db_{wsw^{-1}}=0$ (if $\ell(w)>\ell(ws)$).%
}
\end{itemize}

Here, unlike in Axiom $(A)$, the coefficients depend
only on the reflection $wsw^{-1}$ and on whether
$w$ covers or is covered by $ws$ in the Bruhat poset.

\begin{rem}
The dichotomic dependence on the order between $w$ and $ws$
in the Bruhat poset appears also in the action of the Coxeter generators
on certain well-known fundamental bases.
For the action on the Kazhdan-Lusztig basis of the Hecke algebra see~\cite{KL}.
For the action on the Schubert polynomial basis of the coinvariant algebra
see~\cite[Theorem 3.14 (iii)]{BGG}
and the reformulation in~\cite[Theorem 1]{APR}.
\end{rem}




\begin{thm}\label{t.a1}
Let $(\T,\C)$ be a minimal AY pair for
(the Iwahori-Hecke algebra of) $(W,S)$.
If
$$
a_s(w)=a_{s'}(w') \then b_s(w)=b_{s'}(w')
\qquad(\forall s,s'\in S, w,w'\in \C)
$$
then
$\T$ satisfies Axiom $(B)$.
\end{thm}

This result will be proved in Appendix 1 (Section~\ref{s.ax1}).

\smallskip

%

Thus, under reasonable assumptions, the coefficients
$a_s(w)$ and $b_s(w)$ in Axiom $(A)$ depend only on the reflection
$wsw^{-1}\in T$ and on whether or not $\ell(w)<\ell(ws)$.

\medskip

It is clear that in Axiom $(B)$ we actually use only a subset of
the set $T$ of all reflections.
Define:
\begin{eqnarray*}
T_{\C}          &:=& \{wsw^{-1}\,|\,s\in S,\,w\in \C,\,ws\in \C\},\\
T_{\partial \C} &:=& \{wsw^{-1}\,|\,s\in S,\,w\in \C,\,ws\not\in \C\}.
\end{eqnarray*}

Recall Corollary~\ref{t.convex_gdc}
and the notation $\C_A(w)$ from Definition~\ref{d.A_equiv}.

\begin{lem}\label{t.part}
Let $\C\subseteq W$ be a convex set. Then:
\begin{itemize}
\item[{\bf (i)}]
If $\C=\C_A(w)$ (for $A\subseteq T$ and $w\in W$) then
$T_{\C}\subseteq T\setminus A$ and $T_{\partial \C}\subseteq A$.
In particular,
$$
T_{\C} \cap T_{\partial \C} = \emptyset.
$$
\item[{\bf (ii)}]
For any $t\in T_{\partial \C}$,
the edges labeled $t$ in the Hasse diagram are always encountered
``in the same direction'' (either always up or always down)
when going out of $\C$, namely:
if $(w_1, s_1),(w_2,s_2)\in \C\times S$ satisfy
$$
w_1s_1w_1^{-1} = w_2s_2w_2^{-1} = t
$$
(so that, in particular, $w_1s_1, w_2s_2 \in W\setminus \C$)
then
$\ell(w_1s_1)-\ell(w_1)$ and $\ell(w_2s_2)-\ell(w_2)$ have the same sign.
\end{itemize}
\end{lem}

\noindent{\bf Proof.}

\noindent{\bf (i)}
By Corollary~\ref{t.convex_gdc} there exist $A\subseteq T$ and $w\in W$
such that $\C=\C_A(w)$. If $t\in T_{\C}$ then there exist
$s\in S$ and $u\in\C$ such that $t=usu^{-1}$ and $us\in\C$.
By Proposition~\ref{t.des_diff},
$$
\D_T(u) \bigtriangleup \D_T(us) = \{t\},
$$
where $\bigtriangleup$ denotes symmetric difference.
Since by assumption $\D_A(u)=\D_A(us)=\D_A(w)$, it follows that $t\not\in A$.
On the other hand, if $t\in T_{\partial\C}$ then a similar argument shows that
$t\in A$.

\noindent{\bf (ii)}
Let $\C=\C_A(w)$ as above, and let $t\in T_{\partial\C}\subseteq A$.
Let $(w_1, s_1),(w_2,s_2)\in \C\times S$ satisfy
$w_1s_1w_1^{-1} = w_2s_2w_2^{-1} = t$.
Since $w_1, w_2 \in \C_A(w)$,
$$
\D_A(w_1) = \D_A(w_2) = \D_A(w).
$$
If $t\in \D_A(w_1) = \D_A(w_2)$ then
$\ell(w_i) > \ell(tw_i) = \ell(w_i s_i)$ ($i=1,2$).
If $t\not\in \D_A(w_1) = \D_A(w_2)$ then
$\ell(w_i) < \ell(tw_i) = \ell(w_i s_i)$ ($i=1,2$).
Therefore $\ell(w_1s_1) - \ell(w_1)$ and $\ell(w_2s_2) - \ell(w_2)$
have the same sign.

\qed

\medskip

\begin{df}\label{d.out}
For $t\in T_{\partial \C}$ define $\oa_t$ to be either $\ua_t$ or $\da_t$
in the direction ``out of $\C$'', as in Lemma~\ref{t.part}.
Define $\ob_t$ similarly.
\end{df}

\subsection{Normalization}\label{s.normal}

For the following lemma, we assume only that $\T$ is a map
from the generating set $S$
to the algebra $E_{\C} := End_{\bbf}(V_{\C})$, where $\C$ is a 
subset of $W$,
and that $\T$ satisfies Axiom $(B)$.
Denote $\T_s:=\T(s)$ $(\forall s\in S)$.

\begin{lem}\label{t.l1}
Let $\C$ be a subset of $W$, let
$\T:S \too E_{\C}$ be a map satisfying Axiom $(B)$, and assume
that
$$
\ub_t, \db_t \ne 0\qquad(\forall t\in T_{\C}).
$$
Then the Hecke relation
\begin{equation}\label{e.hecke}
(\T_s-1)(\T_s+q_s)=0\qquad(\forall s\in S)
\end{equation}
holds if and only if
\begin{itemize}
\item[(a)]
For any $t\in T_{\C}$:
$$
\ua_t+\da_t = 1-q_t,\quad
\ub_t \db_t = (1-\ua_t)(1-\da_t) \ne 0.
$$
\item[(b)]
For any $t\in T_{\partial \C}$:
$$
\oa_t\in \{1,-q_t\},\quad
\ob_t = 0.
$$
\end{itemize}
\end{lem}

\noindent
{\bf Proof.} 
It will be more convenient to use in this proof
the notation $a_s(w)$ of Axiom $(A)$ rather than
the notation $\ua_{wsw^{-1}}$ (or $\da_{wsw^{-1}}$) of Axiom $(B)$.

Fix $s\in S$ and $w\in \C$.
If $ws\in \C$ then
\begin{eqnarray*}
(\T_s-1)(\T_s+q_s)(C_w)
&=& (\T_s-1)[(a_s(w)+q_s)C_w + b_s(w)C_{ws}]\\
&=& (a_s(w)+q_s)[(a_s(w)-1)C_w + b_s(w)C_{ws}]\\
&+& b_s(w)[(a_s(ws)-1)C_{ws} + b_s(ws)C_{ws^2}]\\
&=& [(a_s(w)-1)(a_s(w)+q_s)+b_s(w)b_s(ws)]C_w\\
&+& [a_s(w)+a_s(ws)-1+q_s]b_s(w)C_{ws}
\end{eqnarray*}
and, by assumption, $b_s(w)b_s(ws)\ne 0$.
On the other hand, if $ws\not\in \C$ then
$b_s(w)=0$, so that
$$
(\T_s-1)(\T_s+q_s)(C_w)=(a_s(w)-1)(a_s(w)+q_s)C_w.
$$
The deformed involution relation
$$
(\T_s-1)(\T_s+q_s) = 0 \qquad (\forall s\in S)
$$
is thus equivalent to the system of equations
$$
(a_s(w)-1)(a_s(w)+q_s)+b_s(w)b_s(ws)=0\qquad(\forall s\in S, w,ws\in \C),
$$
$$
a_s(w)+a_s(ws)-1+q_s=0\qquad(\forall s\in S, w,ws\in \C),
$$
and
$$
(a_s(w)-1)(a_s(w)+q_s)=0\qquad(\forall s\in S, w\in \C, ws\not\in \C).
$$

Assume first that $ws\in \C$. Since
$$
a_s(w)+a_s(ws) = 1-q_s,
$$
we can replace $a_s(w)+q_s$ by $1-a_s(ws)$ in the first equation to get
$$
(a_s(w)-1)(1-a_s(ws))+b_s(w)b_s(ws) = 0.
$$
If $ws\not\in\C$ then clearly
$$
a_s(w)\in\{1,-q_s\}.
$$
Combining these results with the vanishing properties of $b_s(w)$
completes the proof, upon converting to the notation of Axiom $(B)$
and Definition~\ref{d.out} and recalling (from Section~\ref{s.prelim})
that $q_{wsw^{-1}} = q_s$.
\qed

Note that, if $s\in S$, $w,ws\in\C$, and $\ell(w) < \ell(ws)$,
then $\T_s$ is represented on the invariant subspace
$\sn_{\bbf}\{C_w, C_{ws}\}$ by a $2 \times 2$ matrix as follows:
$$
\left[\begin{array}{c}
\T_s(C_w)\cr \T_s(C_{ws})
\end{array}\right] =
\left[\begin{array}{cc}
\ua_t & \ub_t\cr
\db_t & \da_t
\end{array}\right]
\left[\begin{array}{c}
C_w\cr C_{ws}
\end{array}\right],
$$
where $t:=wsw^{-1}\in T_{\C}$. By Lemma~\ref{t.l1}, this matrix has
trace $1-q_t$ and determinant $-q_t$.
We can further restrict this matrix by adding normalization assumptions,
as follows.

As we shall see in the next section, the coefficients $\ua_t$, $\da_t$
(and $\oa_t$) must satisfy certain constraints, in addition to those of
Lemma~\ref{t.l1}, in order for $\T$ to be a representation.
On the other hand, $\ub_t$ and $\db_t$ have no additional constraints.
It is therefore natural to normalize $\ub_t$ and $\db_t$.
We list a few possibilities.

\begin{description}
\item{$(SON)$}
Symmetric normalization:
$$
\db_t = \ub_t
\qquad(\forall t\in T_{\C}).
$$
Note that if $q_t=1$ then $\da_t = -\ua_t$ and $\ua_t^2 + \ub_t\db_t = 1$,
so that under this normalization the matrix representing
each generator $\T_s$ (in the natural basis of $V_{\C}$) is {\em orthogonal}
as well as symmetric.
Also, due to the quadratic constraint on $\ua_t$ and $\ub_t$,
this normalization is not available for arbitrary fields $\bbf$.
\item{$(SNN)$}
Seminormal normalization:
$$
\db_t = 1
\qquad(\forall t\in T_{\C}).
$$
\item{$(RSN)$}
Row stochastic normalization:
$$
\ua_t + \ub_t = \da_t + \db_t = 1
\qquad(\forall t\in T_{\C}).
$$
\item{$(CSN)$}
Column stochastic normalization:
$$
\ua_t + \db_t = \da_t + \ub_t = 1
\qquad(\forall t\in T_{\C}).
$$
\end{description}

Normalizations $(SON)$ and $(SNN)$,
for the special case of the symmetric group,
were introduced by Alfred Young~\cite{Ja}.
Hoefsmit~\cite{Ho} and Ram~\cite{Ra} used normalization $(CSN)$
in their studies of Hecke algebra representations.

\bigskip

It follows that, in order to define a representation, it suffices to determine the
coefficients $\{\ua_t\,|\,t\in T_{\C}\}$ and $\{\oa_t\,|\,t\in T_{\partial \C}\}$
(and to choose a normalization).
Later we will show that,
for simply laced Coxeter groups, it actually suffices to determine
these coefficients on the subset $\{vsv^{-1}\,|\,s\in S\}\subseteq T$
for some {\em fixed} $v\in \C$; the other coefficients are then determined
by a linearity condition. See Theorems~\ref{t.f2rep} and~\ref{t.rep2f}
below.

\section{Coxeter Relations and Coefficients}\label{s.CRC}

In this section (as in Lemma~\ref{t.l1} above)
we assume only that $\T$ is a map from the generating set $S$
to the algebra $E_{\C} := End_{\bbf}(V_{\C})$,
where $\C$ is a subset (in Lemma~\ref{t.l3}, a convex subset) of $W$,
and that $\T$ satisfies Axiom $(B)$.
Denote $\T_s:=\T(s)$ $(\forall s\in S)$.

We shall find necessary and sufficient conditions for $\T$ to satisfy
the Coxeter relations for $m = 2,3$.  Under these conditions $\T$ will
actually be a representation of ${\cal H}_q(W)$
in case $W$ is simply laced, as we shall see in subsequent sections.

%

\begin{lem}\label{t.l2}
Let $\C$ be a subset of $W$ and let
$\T:S \too E_{\C}$ satisfy Axiom $(B)$.
If $s,t\in S$ satisfy $m(s,t)=2$ (i.e., $st=ts$)
then $\T_s\T_t=\T_t\T_s$.
\end{lem}

\noindent
{\bf Proof.} 
As in the proof of Lemma~\ref{t.l1}, we shall use for convenience
the notation of Axiom $(A)$.
Take $w\in \C$, and assume first that $ws,wt\in\C$. Then:
\begin{eqnarray*}
\T_s \T_t (C_w) &=&
\T_s(a_t(w)C_w + b_t(w)C_{wt}) =\\
&=& a_t(w)a_s(w)C_w
 +  a_t(w)b_s(w)C_{ws} +\\
&+& b_t(w)a_s(wt)C_{wt}
 +  b_t(w)b_s(wt)C_{wts}.
\end{eqnarray*}
Similarly
\begin{eqnarray*}
\T_t\T_s (C_w)
&=& a_s(w)a_t(w)C_w
 +  a_s(w)b_t(w)C_{wt} +\\
&+& b_s(w)a_t(ws)C_{ws}
 +  b_s(w)b_t(ws)C_{wst}.
\end{eqnarray*}
Noting that $wst = wts$, the equation $\T_s\T_t(C_w) = \T_t\T_s(C_w)$
is thus equivalent to the system of equations
$$
\begin{array}{lcrcl}
C_w&:\quad&
a_t(w)a_s(w) &=& a_s(w)a_t(w)\\
C_{ws}&:\quad&
a_t(w)b_s(w) &=& b_s(w)a_t(ws)\\
C_{wt}&:\quad&
b_t(w)a_s(wt) &=& a_s(w)b_t(w)\\
C_{wst}&:\quad&
b_t(w)b_s(wt) &=& b_s(w)b_t(ws)
\end{array}
$$

It is enough to show that
$$
\begin{array}{rcl}
a_t(w) &=& a_t(ws)\\
a_s(wt) &=& a_s(w)\\
b_t(w)b_s(wt) &=& b_s(w)b_t(ws)
\end{array}
$$
The first two equations actually follow from
Axiom $(B)$ together with
the assumption $st=ts$. Indeed, $(ws)t(ws)^{-1} = wstsw^{-1} = wtw^{-1}$
and $ws<wst \iff w<wt$ imply $a_t(ws) = a_t(w)$.
Similarly for $a_s(wt) = a_s(w)$.
The same argument shows that
$b_t(ws) = b_t(w)$ and $b_s(wt) = b_s(w)$, implying the third equation.

Assume now that $w,ws\in\C$ but $wt\not\in\C$.
Then $b_t(w)=0$, and also $b_t(ws)=0$ because of Axiom $(B)$.
Thus $\T_s\T_t(C_w) = \T_t\T_s(C_w)$ is equivalent in this case
to the equations
$$
\begin{array}{lcrcl}
C_w&:\quad&
a_t(w)a_s(w) &=& a_s(w)a_t(w)\\
C_{ws}&:\quad&
a_t(w)b_s(w) &=& b_s(w)a_t(ws)\\
\end{array}
$$
which are clearly satisfied since again $a_t(w)=a_t(ws)$.

Similar arguments apply if $w, wt\in\C$ but $ws\not\in\C$,
or if $w\in\C$ but $ws, wt\not\in\C$.

In all cases, the condition $\T_s\T_t(C_w) = \T_t\T_s(C_w)$ yields
no additional restrictions on the coefficients.

\qed



As a preparation for the next result, note the following.

\begin{obs}\label{t.coset}
Let $w\in W$, and let $s, t\in S$, $s\ne t$.  Then:
\begin{itemize}
\item[(a)]
The coset $w\langle s,t \rangle$ contains a unique element $\tw$ such that
$\ell(\tw) < \ell(\tw s)$ and $\ell(\tw) < \ell(\tw t)$.
$\tw$ is the shortest element in $w\langle s,t \rangle$
(``minimal coset representative'').
\item[(b)]
Assume that $m := m(s,t) < \infty$.
If $\C\subseteq W$ is convex then either $w\langle s,t \rangle \subseteq \C$,
or $\C \cap w\langle s,t \rangle$ is a convex set consisting of
at most $m$ elements.
\end{itemize}
\end{obs}

\begin{lem}\label{t.l3}
Let $\C$ be a convex subset of $W$, let
$\T:S \too E_{\C}$ satisfy Axiom $(B)$ and the Hecke relation~(\ref{e.hecke})
from Lemma~\ref{t.l1}, and assume that
$\ub_t, \db_t \ne 0$ $(\forall t\in T_{\C})$.
If $s,t\in S$ satisfy $m(s,t)=3$ (i.e., $sts=tst$)
then $\T_s\T_t\T_s=\T_t\T_s\T_t$ if and only if:
\begin{itemize}
\item[(a)]
For every $w\in \C$ such that
either $ws\in\C$ or $wt\in\C$ (or both):
\begin{equation}\label{e.l3}
\ua_0 \da_2 = \ua_0 \da_1 + \ua_1 \da_2,
\end{equation}
where
\begin{equation}\label{e.l3.1}
\ua_0 := a_s(w) = a_t(wts),
\end{equation}
\begin{equation}\label{e.l3.2}
\ua_1 := a_t(ws) = a_s(wt),
\end{equation}
\begin{equation}\label{e.l3.3}
\ua_2 := a_s(wst) = a_t(w),
\end{equation}
and
\begin{equation}\label{e.l3.4}
\da_i := 1-q-\ua_i\qquad(i=0,1,2).
\end{equation}
Replacing $w$ by any other element of $\C \cap w\langle s,t \rangle$
gives equivalent equations (\ref{e.l3}).
\item[(b)]
For every $w\in \C$ such that
$ws, wt\not\in \C$:
\begin{equation}\label{e.l3.b}
a_s(w) = a_t(w) \in\{1,-q\}.
\end{equation}
\end{itemize}
\end{lem}

Note that $m(s,t)=3$ implies that $s$ and $t$ are conjugate in $W$, so that
$q_s = q_t$ (denoted here $q$).
Also, the notation $\ua_0$ etc.\ here is not related to
$\ua_t$ from Axiom $(B)$, and in particular does not necessarily imply
$\ell(w) < \ell(ws)$ etc.

\medskip

\noindent
{\bf Proof.} 
Assume first that the full coset $w\langle s,t\rangle$ is contained in $\C$.
Then
\begin{eqnarray*}
\T_s \T_t \T_s(C_w) &=&
    a_s(w)a_t(w)a_s(w)C_w
 +  a_s(w)a_t(w)b_s(w)C_{ws} +\\
&+& a_s(w)b_t(w)a_s(wt)C_{wt}
 +  a_s(w)b_t(w)b_s(wt)C_{wts} +\\
&+& b_s(w)a_t(ws)a_s(ws)C_{ws}
 +  b_s(w)a_t(ws)b_s(ws)C_{w} +\\
&+& b_s(w)b_t(ws)a_s(wst)C_{wst}
 +  b_s(w)b_t(ws)b_s(wst)C_{wsts}.
\end{eqnarray*}
and similarly
\begin{eqnarray*}
\T_t \T_s \T_t(C_w) &=&
    a_t(w)a_s(w)a_t(w)C_w
 +  a_t(w)a_s(w)b_t(w)C_{wt} +\\
&+& a_t(w)b_s(w)a_t(ws)C_{ws}
 +  a_t(w)b_s(w)b_t(ws)C_{wst} +\\
&+& b_t(w)a_s(wt)a_t(wt)C_{wt}
 +  b_t(w)a_s(wt)b_t(wt)C_{w} +\\
&+& b_t(w)b_s(wt)a_t(wts)C_{wts}
 +  b_t(w)b_s(wt)b_t(wts)C_{wtst}.
\end{eqnarray*}
Comparing coefficients we get the equations:
$$
\begin{array}{lcc}
C_w&:\quad&
a_s(w)a_t(w)a_s(w) + b_s(w)a_t(ws)b_s(ws) =\\
&&= a_t(w)a_s(w)a_t(w) + b_t(w)a_s(wt)b_t(wt)\\
C_{ws}&:\quad&
a_s(w)a_t(w)b_s(w) + b_s(w)a_t(ws)a_s(ws) = a_t(w)b_s(w)a_t(ws)\\
C_{wt}&:\quad&
a_s(w)b_t(w)a_s(wt) = a_t(w)a_s(w)b_t(w) + b_t(w)a_s(wt)a_t(wt)\\
C_{wst}&:\quad&
b_s(w)b_t(ws)a_s(wst) = a_t(w)b_s(w)b_t(ws)\\
C_{wts}&:\quad&
a_s(w)b_t(w)b_s(wt) = b_t(w)b_s(wt)a_t(wts)\\
C_{wsts}&:\quad&
b_s(w)b_t(ws)b_s(wst) = b_t(w)b_s(wt)b_t(wts)
\end{array}
$$

Let $w_0 = sts =tst$ be the longest element in the parabolic subgroup
$\langle s,t \rangle$ of $W$. Multiplication (on the right) by $w_0$
in the coset $w\langle s,t \rangle$ is order reversing, so that
$\ell(w)<\ell(ws) \iff \ell(wsw_0)<\ell(ww_0)$. Thus, by Axiom $(B)$,
$$
a_s(w) = a_{w_0sw_0}(wsw_0) = a_t(wts).
$$
This is equality (and notation) (\ref{e.l3.1}) above.
Similarly for (\ref{e.l3.2}), (\ref{e.l3.3}) and also (\ref{e.l3.4})
(since $q_{wsw^{-1}} = q_{wstsw^{-1}} = q_{wtw^{-1}} = q$).

Define similarly $\ub_i$ and $\db_i$ $(0\le i\le 2)$.
Using this notation, we get
$$
\begin{array}{lcrcl}
C_w&:\quad&
\ua_0 \ua_2 \ua_0 + \ub_0 \ua_1 \db_0 &=&
\ua_2 \ua_0 \ua_2 + \ub_2 \ua_1 \db_2\\
C_{ws}&:\quad&
\ua_0 \ua_2 \ub_0 + \ub_0 \ua_1 \da_0 &=& \ua_2 \ub_0 \ua_1\\
C_{wt}&:\quad&
\ua_0 \ub_2 \ua_1 &=& \ua_2 \ua_0 \ub_2 + \ub_2 \ua_1 \da_2\\
C_{wst}&:\quad&
\ub_0 \ub_1 \ua_2 &=& \ua_2 \ub_0 \ub_1\\
C_{wts}&:\quad&
\ua_0 \ub_2 \ub_1 &=& \ub_2 \ub_1 \ua_0\\
C_{wsts}&:\quad&
\ub_0 \ub_1 \ub_2 &=& \ub_2 \ub_1 \ub_0
\end{array}
$$
The last three equations are tautologies. The first three become
(after division by $1$, $\ub_0 \ne 0$ and $\ub_2 \ne 0$, respectively):
$$
\begin{array}{lcrcl}
C_w&:\quad&
\ua_0 \ua_2 (\ua_0 - \ua_2) &=& \ua_1 (\ub_2 \db_2 - \ub_0 \db_0)\\
C_{ws}&:\quad&
\ua_0 \ua_2 + \da_0 \ua_1 &=& \ua_1 \ua_2\\
C_{wt}&:\quad&
\ua_0 \ua_1 &=& \ua_0 \ua_2 + \ua_1 \da_2\\
\end{array}
$$
From equation~(\ref{e.l3.4})
it follows that
$$
\ua_1 (\ua_0 + \da_0) = \ua_1 (\ua_2 + \da_2).
$$
Subtract this from the equation for $C_{ws}$ to get, equivalently,
$$
\ua_0 \ua_2 - \ua_0 \ua_1 = -\ua_1 \da_2.
$$
Now use similarly
$$
\ua_0 (\ua_2 + \da_2) = \ua_0 (\ua_1 + \da_1)
$$
to get the claimed equation~(\ref{e.l3}).

Similar operations on the equation for $C_{wt}$ lead to
$$
-\ua_0 \da_1 = -\ua_0 \da_2 + \ua_1 \da_2,
$$
which is again equation~(\ref{e.l3}).
Finally, from Lemma~\ref{t.l1} and equation~(\ref{e.l3.4}):
$$
\ub_2 \db_2 - \ub_0 \db_0 = (1-\ua_2)(1-\da_2) - (1-\ua_0)(1-\da_0) =
\ua_2 \da_2 - \ua_0 \da_0 = (\ua_0 - \ua_2)(\ua_2 - \da_0).
$$
Thus the equation for $C_w$ is equivalent to
$$
\ua_0 \ua_2 (\ua_0 - \ua_2) = \ua_1 (\ua_0 - \ua_2)(\ua_2 - \da_0)
$$
which is clearly a consequence of
$$
\ua_0 \ua_2 = \ua_1 (\ua_2 - \da_0),
$$
again equivalent to equation~(\ref{e.l3}).

Assume now that $w\langle s,t\rangle\not \subseteq \C$.
By Observation~\ref{t.coset}(b), $\C\cap w\langle s,t\rangle$ is
a convex set containing at most 3 elements. In particular, $wsts\not\in\C$.

If $\C \cap w\langle s,t \rangle = \{w,ws,wst\}$ then $\ub_2=0$,
and equations $C_w$, $C_{ws}$ and $C_{wst}$ above are
$$
\begin{array}{lcrcl}
C_w&:\quad&
\ua_0 \ua_2 \ua_0 + \ub_0 \ua_1 \db_0 &=& \ua_2 \ua_0 \ua_2\\
C_{ws}&:\quad&
\ua_0 \ua_2 \ub_0 + \ub_0 \ua_1 \da_0 &=& \ua_2 \ub_0 \ua_1\\
C_{wst}&:\quad&
\ub_0 \ub_1 \ua_2 &=& \ua_2 \ub_0 \ub_1\\
\end{array}
$$
Again, the third equation is trivial; the second is equivalent to
equation~(\ref{e.l3}) (since $\ub_0 \ne 0$); and the first is a consequence
of equation~(\ref{e.l3}).

The cases $\C \cap w\langle s,t \rangle = \{w,wt,wts\}$ (with $\ub_0=0$)
and  $\C \cap w\langle s,t \rangle = \{w,ws,wt\}$ (with $\ub_1=0$)
are similar.

If $\C \cap w\langle s,t \rangle = \{w,ws\}$ then $\ub_1=\ub_2=0$.
Equations $C_w$ and $C_{ws}$ are
$$
\begin{array}{lcrcl}
C_w&:\quad&
\ua_0 \ua_2 \ua_0 + \ub_0 \ua_1 \db_0 &=& \ua_2 \ua_0 \ua_2\\
C_{ws}&:\quad&
\ua_0 \ua_2 \ub_0 + \ub_0 \ua_1 \da_0 &=& \ua_2 \ub_0 \ua_1\\
\end{array}
$$
Again, the second equation is equivalent to
equation~(\ref{e.l3}) (since $\ub_0 \ne 0$), and the first is a consequence
of it.

The case $\C \cap w\langle s,t \rangle = \{w,wt\}$ is similar.

To complete the proof of part (a) of the lemma we must show that
replacing $w$ by any other element of $\C \cap w\langle s,t \rangle$
yields an equivalent equation~(\ref{e.l3}).
Due to the convexity of $\C \cap w\langle s,t \rangle$, this will follow by
induction once we prove it for the replacement of $w$ by $ws$ or by $wt$.
Assume, e.g., that $w, ws\in \C$. Replacing $w$ by $ws$ gives
$$
\begin{array}{lllll}
\hbox{\rm new\ } \ua_0 &=& a_s(ws) &=& \hbox{\rm old\ } \da_0,\cr
\hbox{\rm new\ } \ua_1 &=& a_t(w) &=& \hbox{\rm old\ } \ua_2,\cr
\hbox{\rm new\ } \ua_2 &=& a_s(wt) &=& \hbox{\rm old\ } \ua_1.\cr
\end{array}
$$
Thus
$$
\ua_0 \da_2 = \ua_0 \da_1 + \ua_1 \da_2\qquad\hbox{\rm(new)}
$$
is actually
$$
\da_0 \da_1 = \da_0 \da_2 + \ua_2 \da_1\qquad\hbox{\rm(old)},
$$
and by operations as above this can be shown to be equivalent to
the (old) equation~(\ref{e.l3}).

Finally, for part (b):
if $\C \cap w\langle s,t \rangle = \{w\}$ then $\ub_0=\ub_2=0$,
and $\ub_1$ is actually undefined. Equation $C_w$ is now
$$
\begin{array}{lcrcl}
\ua_0 \ua_2 \ua_0 &=& \ua_2 \ua_0 \ua_2\\
\end{array}
$$
Recall that $\ub_0=0$ implies $\ua_0\in\{1,-q\}$, and similarly for $\ua_2$.
In particular $\ua_0,\ua_2\ne 0$, and we conclude from $C_w$ that
$$
\ua_0 = \ua_2,
$$
as claimed in part (b) of the lemma.

\qed

%




\section{Simply Laced Coxeter Groups}\label{s.sl}

In this section, $(W,S)$ is an irreducible simply laced Coxeter system.
Thus its Dynkin diagram is connected, and contains only
``simple'' edges (corresponding to $m(s,t) = 3$) and
non-edges (corresponding to $m(s,t) = 2$).

\medskip

We shall see (Theorem~\ref{t.rep2f}) that, for $q=1$, representations
satisfying Axiom $(B)$ are determined by linear functionals
on the root space.
The $q$-analogues are more subtle; see Section~\ref{s.sl_q}.

\begin{rem}\label{r.unitcell}
By Observation~\ref{t.unitcell}, if $\C$ is a minimal AY cell then,
for any $v\in W$, the left translate $v^{-1}\C$ is a minimal AY cell carrying
an equivalent AY representation (with exactly the same representation matrices).
From now on (unless otherwise stated) we shall assume,
with no loss of generality, that $\C$ contains the identity element of $W$.
This implies, in particular, that
$$
\oa_t = \ua_t \qquad(\forall t\in T_{\partial\C})
$$
since the direction ``out of $\C$'' is always ``upwards'', by convexity.
\end{rem}

We first restate Lemma~\ref{t.l3} in another form, assuming
$\ua_t \ne 0$ ($\forall t\in T_{\C}$).

\begin{lem}\label{t.l3q1}
Let $q=1$, let $\C$ be a convex subset of $W$ containing the identity element,
let $\T:S \to E_{\C}$ satisfy Axiom $(B)$ and the involution relation
$$
\T_s^2 = 1\qquad(\forall s\in S),
$$
and assume that
$$
\ua_t, \ub_t, \db_t \ne 0\qquad(\forall t\in T_{\C}).
$$
If $s,t\in T$ satisfy $m(s,t)=3$,
then $\T_s\T_t\T_s = \T_t\T_s\T_t$ if and only if:
\begin{itemize}
\item[(a)]
For every $w\in\C$ such that $\ell(w)<\ell(ws)$, $\ell(w)<\ell(wt)$,
and either $ws\in\C$ or $wt\in\C$ (or both):
$$
{1\over \ua_{wstsw^{-1}}} =
{1\over \ua_{wsw^{-1}}} + {1\over \ua_{wtw^{-1}}}.
$$
\item[(b)]
For every $w\in\C$ such that $ws,wt\not\in\C$
(and therefore $\ell(w) < \ell(ws)$ and $\ell(w) < \ell(wt)$):
$$
\ua_{wsw^{-1}} = \ua_{wtw^{-1}} = \pm 1.
$$
\end{itemize}
\end{lem}

\noindent
{\bf Proof.}
First note that, since $\C$ is convex and contains the identity element,
if $w\in \C$ then $\C$ also contains the shortest element $\tw$ of
the coset $w\langle s,t \rangle$ (see Observation~\ref{t.coset}(a)).
Since equation~(\ref{e.l3}) of Lemma~\ref{t.l3} is independent of the choice
of element in $\C \cap w\langle s,t \rangle$, we may as well choose $w = \tw$,
namely assume that $\ell(w) < \ell(ws)$ and $\ell(w) < \ell(wt)$.
Lemma~\ref{t.l1} is applicable (since $\ub_t, \db_t \ne 0$ for $t\in T_{\C}$)
and gives, for $q=1$:
$$
\da_t = -\ua_t \qquad(\forall t\in T_{\C})
$$
and
$$
\ua_t = \pm1 \qquad(\forall t\in T_{\partial\C}).
$$
Substitute these into Lemma~\ref{t.l3}, and divide the simplified equation
$$
\ua_0 \ua_2 = \ua_0 \ua_1 + \ua_1 \ua_2
$$
by $\ua_0 \ua_1 \ua_2 \ne 0$.

\qed

Fix an arbitrary positive definite bilinear form $\langle\cdot,\cdot\rangle$
on the root space $V$.

\begin{df}\label{d.cg} ($\C$-genericity)\\
Let $\C$ be a convex subset of $W$ containing the identity element.
A vector $f\in V$ is {\rm $\C$-generic} if:
\begin{itemize}
\item[{\bf (i)}]
For all $t\in T_{\C}$,
$$
\langle f,\al_t\rangle \not\in \{0, 1, -1\}.
$$
\item[{\bf (ii)}]
For all $t\in T_{\partial \C}$,
$$
\langle f,\al_t\rangle = \pm 1.
$$
\item[{\bf (iii)}]
If $w\in \C$, $s,t\in S$, $m(s,t)=3$ and $ws, wt\not\in \C$ then
$$
\langle f,\al_{wsw^{-1}}\rangle =
\langle f,\al_{wtw^{-1}}\rangle\;(= \pm1).
$$
\end{itemize}
\end{df}


The following two theorems hold for all (not necessarily finite) irreducible
simply laced Coxeter groups.

\begin{thm}\label{t.f2rep}
Let $W$ be an irreducible simply laced Coxeter group
and let $\C$ be a convex subset of $W$ containing the identity element.
If $f\in V$ is $\C$-generic then
$$
\ua_t := {1\over \langle f,\al_t\rangle}\qquad
(\forall t\in T_{\C} \cup T_{\partial \C}),
$$
together with $\da_t$, $\ub_t$ and $\db_t$ satisfying
conditions (a) and (b) of Lemma~\ref{t.l1} (with $q = 1$),
define a representation $\T$ of $W$ such that
$(\T,\C)$ is a minimal AY pair satisfying Axiom $(B)$.
\end{thm}

\begin{thm}\label{t.rep2f}
Let $W$ be an irreducible simply laced Coxeter group
and let $\C$ be a subset of $W$ containing the identity element.
If $(\T,\C)$ is a minimal AY pair satisfying Axiom $(B)$
and also $\ua_t\ne 0$ $(\forall t\in T_{\C})$,
then there exists a $\C$-generic $f\in V$ such that
$$
\ua_t={1\over \langle f,\al_t\rangle}\qquad
(\forall t\in T_{\C} \cup T_{\partial \C}).
$$
\end{thm}

\bigskip

\noindent{\bf Proof of Theorem~\ref{t.f2rep}.}
Given a $\C$-generic $f\in V$ define
$$
\ua_t := \frac{1}{\langle f,\al_t\rangle}\qquad
(\forall t\in T_{\C} \cup T_{\partial \C})
$$
and
$$
\da_t := -\ua_t \qquad(\forall t\in T_{\C}).
$$
Define $\ub_t, \db_t\in\bbf$ such that
$$
\ub_t \db_t = (1 - \ua_t)(1 - \da_t)\qquad(\forall t\in T_{\C})
$$
and
$$
\ub_t = 0\qquad(\forall t\in T_{\partial \C}).
$$
Define $\T_s$ (for $s\in S$) by Axiom $(B)$.
It suffices to show that all the defining relations of $W$ are satisfied
by $\{\T_s\,|\,s\in S\}$.
This will be done using Lemmas \ref{t.l1}, \ref{t.l2} and \ref{t.l3q1}.

First of all, $\ua_t$ are defined and $\ua_t, \da_t \ne \pm1$
$(\forall t\in T_{\C})$ because of condition (i) in Definition~\ref{d.cg}.
This implies, by our construction, that $\ub_t, \db_t \ne 0$
$(\forall t\in T_{\C})$ and therefore Lemma~\ref{t.l1} applies.
$\ua_t = \pm1$ $(\forall t\in T_{\partial\C})$ because of
Definition~\ref{d.cg}(ii).

Lemma~\ref{t.l2} clearly applies.

Condition $(a)$ of Lemma~\ref{t.l3q1} is equivalent to
$$
\langle f,\al_{wstsw^{-1}} \rangle =
\langle f,\al_{wsw^{-1}} \rangle +
\langle f,\al_{wtw^{-1}} \rangle
$$
for all $w\in \C$ and $s,t\in S$ such that $m(s,t)=3$, $\ell(w)<\ell(ws)$,
$\ell(w)<\ell(wt)$,
and $|\C\cap w\langle s,t\rangle|\ge 2$.
But for such $w$, $s$ and $t$
$$
\al_{sts}=\sigma_s(\al_t)=\al_t-2B(\al_t,\al_s)\al_s = \al_t+\al_s,
$$
so that
$$
\sigma_w(\al_{sts})=\sigma_w(\al_t+\al_s)=\sigma_w(\al_t)+\sigma_w(\al_s)
$$
or, in other words (since also $\ell(ws) < \ell(wst)$),
$$
\al_{wstsw^{-1}}=\al_{wtw^{-1}}+\al_{wsw^{-1}}
$$
and we are done.

Finally, if $w\in \C$ but $ws,wt\not\in \C$ then
$wsw^{-1}, wtw^{-1}\in T_{\partial \C}$. By Definition~\ref{d.cg}(iii),
condition $(b)$ of Lemma~\ref{t.l3q1} is satisfied.
Minimality of $(\T,\C)$ follows from Corollary~\ref{t.conn_min}.
This completes the proof.

\qed

\bigskip

\noindent{\bf Proof of Theorem~\ref{t.rep2f}.}
Since $\{\al_s\,|\, s\in S\}$ is a basis for $V$,
there exists a unique vector $f\in V$ such that
$$
\langle f, \al_s \rangle = \frac{1}{\ua_s}
\qquad (\forall s \in S).
$$
We claim that
$$
\langle f,\al_t \rangle = \frac{1}{\ua_t}
\qquad (\forall t\in T_{\C} \cup T_{\partial \C})
$$
or, equivalently,
$$
\langle f,\al_{wsw^{-1}} \rangle = \frac{1}{\ua_{wsw^{-1}}}
\qquad (\forall w\in \C, s\in S).
$$
This will be proved by induction on the length of $w$.
Clearly, the claim holds if this length is zero (i.e., if $w=id$).
Suppose that $w \neq id$. We may assume that $\ell(w)<\ell(ws)$,
since otherwise we can replace $w$ by $ws$
(and $id,w\in\C \then ws\in\C$ by convexity).
Let $s_1, \ldots ,s_k \in S$ be such that $s_1 \cdots s_k$ is a reduced word
for $w$. Let $u := s_{1}\ldots s_{k-1}$.
Note that $s\ne s_k$, by the assumption $\ell(w)<\ell(ws)$.
Distinguish two cases.

\medskip

\begin{description}
\item[a)] $m(s,s_k)=2$.
\end{description}
In this case
$$
wsw^{-1} = us_kss_ku^{-1} = usu^{-1}.
$$
Since $u \in \C$ (by convexity of $\C$), it follows by our
induction hypothesis that
$$
\ua_{wsw^{-1}} = \ua_{usu^{-1}}
= \frac{1}{\langle f, \al_{usu^{-1}} \rangle}
= \frac{1}{\langle f, \al_{wsw^{-1}} \rangle},
$$
as desired.

\medskip

\begin{description}
\item[b)] $m(s,s_{k})=3$.
\end{description}

Note that, by definition, $\ell(u)<\ell(us_k)$.
If $\ell(u)>\ell(us)$ then, by convexity, $us\in\C$.
Denoting $v:=us$ we have
$$
wsw^{-1} = us_kss_ku^{-1} = uss_ksu^{-1} = vs_kv^{-1}
$$
so that, by the induction hypothesis (since $\ell(v)<\ell(w)$):
$$
\ua_{wsw^{-1}} = \ua_{vs_kv^{-1}}
= \frac{1}{\langle f,\al_{vs_kv^{-1}} \rangle}
= \frac{1}{\langle f,\al_{wsw^{-1}} \rangle}.
$$
It remains to consider the case $\ell(u)<\ell(us)$.
In this case, by Lemma~\ref{t.l3q1},
$$
\frac{1}{\ua_{wsw^{-1}}} = \frac{1}{\ua_{us_kss_ku^{-1}}}
= \frac{1}{\ua_{us_ku^{-1}}} + \frac{1}{\ua_{usu^{-1}}}.
$$
Since $\ell(u)<\ell(w)$ we have, by the induction hypothesis,
$$
\frac{1}{\ua_{us_{k}u^{-1}}} = \langle f, \al_{us_{k}u^{-1}} \rangle
$$
and
$$
\frac{1}{\ua_{usu^{-1}}} = \langle f, \al_{usu^{-1}} \rangle.
$$
Furthermore, $u$ is the shortest element in the coset $u\langle s,s_k \rangle$
of $W$, and therefore $\ell(us_k)<\ell(us_ks)$. Thus
$$
\al_{wsw^{-1}} = \al_{us_kss_ku^{-1}} = \sigma_{us_k}(\al_s)
= \sigma_u \sigma_{s_k}(\al_s) =
$$
$$
= \sigma_u(\al_s - 2B(\al_s,\al_{s_k})\al_{s_k})
= \sigma_u(\al_{s}+\al_{s_k}) = \al_{usu^{-1}} + \al_{us_ku^{-1}}.
$$
Hence
$$
\frac{1}{\ua_{wsw^{-1}}} =
\frac{1}{\ua_{us_ku^{-1}}} + \frac{1}{\ua_{usu^{-1}}} =
\langle f, \alpha _{us_{k}u^{-1}} \rangle   +
\langle f, \alpha _{usu^{-1}} \rangle  = \langle f, \alpha _{wsw^{-1}} \rangle  ,
$$
as desired.

This concludes the induction step.
Finally, $\C$-genericity of $f$ follows from Lemmas~\ref{t.l1} and~\ref{t.l3q1}.

\qed


\section{Simply Laced Hecke Algebras}\label{s.sl_q}

In this section we state $q$-analogues of Theorems~\ref{t.f2rep}
and \ref{t.rep2f}.
Let $(\T,\C)$ be a minimal AY pair for ${\cal H}_q(W)$ satisfying
Axiom $(B)$.

\medskip

Recall that, for an irreducible simply laced Coxeter system $(W,S)$,
a single value $q_s=q$ ($\forall s\in S$) is used in the definition of
the Hecke algebra ${\cal H}_q(W)$. The field $\bbf$ can be, for example,
the field of rational functions $\bbc(q)$.


We first need a $q$-analogue of Lemma~\ref{t.l3q1}. An apparent obstacle is
that, since $\da_t = (1-q) - \ua_t$, equation (\ref{e.l3}) of Lemma~\ref{t.l3}
now contains (linear and) quadratic terms.

\begin{df}\label{d.q-d}
For any reflection $t\in T_{\C} \cup T_{\partial \C}$ such that $\ua_t\ne 0$
let
$$
\ud_t := 1 - \frac{1-q}{\ua_t}.
$$
\end{df}
Thus
$$
\ua_t = \frac{1-q}{1-\ud_t}.
$$

\begin{pro}\label{t.q-prop}
If $q\ne 1$ and
$\ua_0, \ua_1, \ua_2 \ne 0$ (in the notation of Lemma~\ref{t.l3})
then equation~(\ref{e.l3}) of Lemma~\ref{t.l3} is equivalent to
$$
\ud_1 = \ud_0 \cdot \ud_2.
$$
\end{pro}

\medskip

\noindent{\bf Proof.}
We assume that $q\ne 1$ and
$\ua_r\ne 0$ for $r\in\{0,1,2\}$.
thus $\ud_r$ is well-defined and $\ud_r\ne 1$ (for these values of $r$).
If we define similarly
$$
\dd_r := 1 - \frac{1-q}{\da_r}
$$
then
$
\ua_r + \da_r = 1-q
$
is equivalent to
$$
\frac{1}{1-\ud_r} + \frac{1}{1-\dd_r} = 1,
$$
namely to
$$
\ud_r \dd_r = 1.
$$
Equation~(\ref{e.l3}) of Lemma~\ref{t.l3} is now equivalent to
$$
\frac{1}{1-\ud_0} \cdot \frac{1}{1-\dd_2} =
\frac{1}{1-\ud_0} \cdot \frac{1}{1-\dd_1} +
\frac{1}{1-\ud_1} \cdot \frac{1}{1-\dd_2},
$$
namely to
$$
\frac{-\ud_2}{(1-\ud_0)(1-\ud_2)} =
\frac{-\ud_1}{(1-\ud_0)(1-\ud_1)} +
\frac{-\ud_2}{(1-\ud_1)(1-\ud_2)}.
$$
Clearing denominators gives
$$
\ud_2(1-\ud_1) = \ud_1(1-\ud_2) + \ud_2(1-\ud_0)
$$
which is the claimed equation
$$
\ud_1 = \ud_0 \cdot \ud_2.
$$

\qed

Thus, the proper $q$-analogue of the additive condition
on $\frac{1}{\ua_r}$ for $q=1$ (as in Lemma~\ref{t.l3q1})
is a multiplicative condition on $\ud_r$.

\begin{lem}\label{t.l3q} {\rm ($q$-analogue of Lemma~\ref{t.l3q1})}
Let $\C$ be a convex subset of $W$ containing the identity element,
let $\T:S \to E_{\C}$ satisfy Axiom $(B)$ and the Hecke relation
$$
(\T_s - 1)(\T_s + q) = 0\qquad(\forall s\in S),
$$
and assume that
$$
\ua_t, \ub_t, \db_t \ne 0\qquad(\forall t\in T_{\C}).
$$
Let $\ud_t$ be as in Definition~\ref{d.q-d}.
If $s,t\in T$ satisfy $m(s,t)=3$,
then $\T_s\T_t\T_s = \T_t\T_s\T_t$ if and only if:
\begin{itemize}
\item[(a)]
For every $w\in\C$ such that $\ell(w)<\ell(ws)$, $\ell(w)<\ell(wt)$,
and either $ws\in\C$ or $wt\in\C$ (or both):
$$
\ud_{wstsw^{-1}} = \ud_{wsw^{-1}} \cdot \ud_{wtw^{-1}}.
$$
\item[(b)]
For every $w\in\C$ such that $ws,wt\not\in\C$
(and therefore $\ell(w) < \ell(ws)$ and $\ell(w) < \ell(wt)$):
$$
\ud_{wsw^{-1}} = \ud_{wtw^{-1}} \in \{q^1,q^{-1}\}.
$$
\end{itemize}
\end{lem}
The proof is similar to that of Lemma~\ref{t.l3q1}, and will be omitted.

\bigskip

\noindent
For any integer $k\in\bbz$ and parameter $q$ let
$$
[k]_q:=\cases{%
\frac{1-q^k}{1-q},& if $q\ne 1$;\cr
k,& if $q=1$.}
$$
For any $k\in\bbz$, $[k]_q$ is a rational function (actually, a Laurent
polynomial) of $q$. In particular: $[-1]_q=-q^{-1}$, $[1]_q=1$ and $[0]_q=0$.

\smallskip



\begin{thm}\label{t.f2rep.q} {\rm ($q$-analogue of Theorem~\ref{t.f2rep})}
Let ${\cal H}_q(W)$ be the Iwahori-Hecke algebra of
an irreducible simply laced Coxeter group $W$,
and let $\C$ be a convex subset of $W$ containing the identity element.
If $f\in V$ is $\C$-generic
and $\langle f,\al_s\rangle \in\bbz$ $(\forall s\in S)$,
then
$$
\ua_t := {1\over [\langle f,\al_t\rangle]_q}\qquad
(\forall t\in T_{\C} \cup T_{\partial \C}),
$$
together with $\da_t$, $\ub_t$ and $\db_t$ satisfying
conditions (a) and (b) of Lemma~\ref{t.l1},
define a representation $\T$ of ${\cal H}_q(W)$ such that
$(\T,\C)$ is a minimal AY pair satisfying Axiom $(B)$.
\end{thm}

\noindent{\bf Proof.}
Analogous to the proof of Theorem~\ref{t.f2rep}.
Since $\ua_r={1\over [\langle f, \al_r\rangle]_q}$,
$\ud_r=q^{\langle f, \al_r\rangle}$.
Condition $(a)$ of Lemma~\ref{t.l3} is satisfied, by
Proposition~\ref{t.q-prop},
since
$$
\ud_{wstsw^{-1}} = q^{\langle f, \al_{wstsw^{-1}}\rangle} =
q^{\langle f, \al_{wsw^{-1}} + \al_{wtw^{-1}}\rangle} =
\ud_{wsw^{-1}} \cdot \ud_{wtw^{-1}}.
$$

\qed


\bigskip


\begin{thm}\label{t.rep2f.q} {\rm ($q$-analogue of Theorem~\ref{t.rep2f})}
Let ${\cal H}_q(W)$ be the Iwahori-Hecke algebra of
an irreducible simply laced Coxeter group $W$,
and let $\C$ be a subset of $W$ containing the identity element.
If $(\T,\C)$ is a minimal AY pair for ${\cal H}_q(W)$ satisfying Axiom $(B)$
such that $\ua_t\ne 0$ $(\forall t\in T_{\C})$
and $\ud_s$ (as in Definition~\ref{d.q-d}) is an integral power of $q$
$(\forall s\in S)$,
then there exists a $\C$-generic $f\in V$ such that
$$
\ua_t= \frac{1}{[\langle f,\al_t\rangle]_q}\qquad
(\forall t\in T_{\C} \cup T_{\partial \C}).
$$
\end{thm}

\noindent{\bf Proof.}
By induction on length as for $q=1$, replacing the additivity of $\frac{1}{\ua_r}$
by the multiplicativity of $\ud_r$.

\qed

Note that Theorems~\ref{t.f2rep.q} and~\ref{t.rep2f.q} have more general
versions (with less restrictive assumptions), but their formulations are
more complicated.

%

\section{Induction and Restriction}\label{s.induce}

The following observation gives a combinatorial procedure for
restricting AY representations to parabolic subgroups. This
procedure is analogous to the one given by Barbasch and Vogan for
Kazhdan-Lusztig cell representations~\cite[Prop. 3.11]{BV2}.

\begin{obs}\label{t.i1}
Let $(W,S)$ be a finite Coxeter system, and let $P=\langle J\rangle$
$(J\subseteq S)$ be a parabolic subgroup of $W$.
Let $(\rho,\C)$ be a minimal AY pair for $W$. Then
\begin{itemize}
\item[(1)]
$\C$ is a disjoint union of sets of the form $r_i{\cal D}_i$ where,
for each $i$, $r_i\in W$ and ${\cal D}_i$ is a minimal AY cell for $P$.
\item[(2)]
The restricted representation ${\rho\downarrow}_P^W$ is isomorphic to
the direct sum $\bigoplus_{i}\psi^{{\cal D}_i}$, where the
sum runs over the ${\cal D}_i$ in (1) above and, for each $i$,
$\psi^{{\cal D}_i}$ is a minimal AY representation of $P$ on $V_{{\cal D}_i}$.
\end{itemize}
\end{obs}

\noindent{\bf Proof.}
Axiom $(A)$ is equivalent to the following statement:
for any parabolic subgroup $P = \langle J\rangle$ of $W$
(where $J\subseteq S$) and any $w\in \C$, the ``coset subspace''
$V_{wP \cap \C}$ is invariant under $\T(p)$, for all $p\in P$
(see Observation~\ref{t.aprime} above).
\qed

\begin{rem}
Here $\psi^{{\cal D}_i}$ is some AY representation on the AY cell
${\cal D}_i$. Note that, unlike the analogous restriction rule for
Kazhdan-Lusztig representations, $\psi^{{\cal D}_i}$
is not uniquely determined by ${\cal D}_i$.
\end{rem}

The following is a combinatorial procedure for induction, which is analogous to
the one for Kazhdan-Lusztig cells,
proved by Barbasch and Vogan for Weyl groups~\cite[Prop. 3.15]{BV2}
and by Geck for Coxeter groups~\cite{Geck}.

\begin{thm}\label{t.i4}
Let $(W,S)$ be a finite Coxeter system, $P=\langle J\rangle$
$(J\subseteq S)$ a parabolic subgroup of $W$,
and $W^J$ be the set of all representatives of minimal length of
the right cosets of $P$ in $W$.
If $(\psi,\cal D)$ is a minimal AY pair for $P$, then:
\begin{itemize}
\item[1.]
${\cal D} W^J$ is a minimal AY cell for $W$.
\item[2.]
The induced representation $\psi\uparrow ^W_P$ is isomorphic to
a (minimal) AY representation on $V_{{\cal D} W^J}$.
\end{itemize}
\end{thm}

\begin{rem}
${\cal D} W^J$ is a minimal AY cell, and not a just a union of AY
cells as in the analogous theorem for Kazhdan-Lusztig cells.
\end{rem}

For the proof of Theorem~\ref{t.i4} we shall need the following

%

\begin{pro}{\rm\cite[Proposition~3.1.3]{BB}}
For $r_1,r_2\in W$, $r_1\le r_2$ (in weak Bruhat order)
if and only if $\D_T(r_1)\subseteq \D_T(r_2)$.
\end{pro}

Since $r$ belongs to $W^J$ if and only if $\D_T(r)\cap J=\emptyset$,
we conclude

\begin{cor}\label{t.fact}
If $r\in W^J$ and $r'\le r$ in weak Bruhat order,
then $r'\in W^J$.
\end{cor}

%

\begin{lem}\label{t.i5}
Let $(W,S)$ be a finite Coxeter system, $P$ a parabolic subgroup
of $W$ generated by $J\subseteq S$ and $W^J$ the set of minimal
right coset representatives of $P$ in $W$. Let $s\in S$ be a
simple reflection and let $r\in W^J$. Then: either
\begin{itemize}
\item[(1)]
$rs\in W^J$; or
\item[(2)]
$rs\not\in W^J$ and $rs=pr$,
where $p\in J = S\cap P$ is a simple reflection in $P$.
\end{itemize}
\end{lem}

\noindent{\bf Proof.}
If $\ell(rs) < \ell(r)$ then
$\D_T(rs) = \D_T(r)\setminus\{rsr^{-1}\} \subseteq \D_T(r)$
and therefore $rs\in W^J$.
Assume, therefore, that $\ell(rs) > \ell(r)$. Then
$\D_T(rs) = \D_T(r)\cup\{rsr^{-1}\}$.
If $rsr^{-1}\not\in J$ then $\D_T(rs) \cap J = \emptyset$,
so that again $rs\in W^J$.
If $rsr^{-1}\in J$ then $rs=pr$ with $p=rsr^{-1}\in J$,
as claimed.
\qed

\medskip

\noindent{\bf Proof of Theorem~\ref{t.i4}.}
Let $V_{\cal D}=\sn_{\bbc}\{C_m\,|\,m\in {\cal D}\}$
be the representation space of $\psi$.
$V_{\cal D}$ is a right $\bbc[P]$-module.
A representtaion space for $\rho=\psi\uparrow_P^W$ is
the right $\bbc[W]$-module $N=V_{\cal D}\otimes_{\bbc[P]} \bbc[W]$,
with basis $B=\{C_m\otimes r\,|\,m\in {\cal D}, r\in W^J\}$.
The action of $\rho$ on basis elements is
$$
\rho_s(C_m\otimes r)=C_m\otimes rs.
$$
There are two cases:
\begin{itemize}
\item[(1)]
$rs\in W^J$: Then $\rho_s(C_m\otimes r)=C_m\otimes rs\in B$.
\item[(2)]
$rs\not\in W^J$: Then, by Lemma~\ref{t.i5},
$rs=pr$ with $p\in J=P\cap S$. Since $({\cal D},\psi)$ is an AY pair
for $P$, we can write
$$
\psi_p(C_m) = a_p(m) C_m + b_p(m) C_{mp}.
$$
Thus
\begin{eqnarray*}
\rho_s(C_m \otimes r)
&=& C_m \otimes rs = C_m \otimes pr =\cr
&=& \psi_p(C_m) \otimes r = a_p(m) C_m \otimes r + b_p(m) C_{mp} \otimes r.
\end{eqnarray*}
\end{itemize}
Using the natural bijection
$\phi: B\to \{C_{mr}\,|\,m\in {\cal D}, r\in W^J\}$
given by $\phi(C_m\otimes r) := C_{mr}$ we now have:
$$
\rho_s(C_{mr}) = \cases{%
C_{mrs}, &if $rs\in W^J$;\cr
a_p(m)C_{mr} + b_p(m)C_{mrs}, &otherwise.}
$$
We conclude that the induced representation $\rho$ on $V_{{\cal
D}W^J}=\sn\{C_{mr} |\ m\in {\cal D}, r\in W^J\}$ is AY. Thus
${\cal D}W^J$ is an AY cell in $W$.

Verifying minimality (i.e., strong connectivity) of this cell
is easy and is left to the reader. \qed

\noindent{\bf Historical Note:}
Our first proof of Theorem~\ref{t.i4} relied on the fact that, for finite $W$,
$W^J$ is a lattice under weak Bruhat order. This was proved by Bj\"orner and
Wachs~\cite[Theorem 4.1]{BW}.

\section{Examples}\label{examples}

Several concrete examples of AY representations will be described
in this section.

\subsection{Descent Representations}\label{s.des_rep}

In this subsection we consider standard (rather than generalized)
descent classes, and assume $q=1$. It is possible to carry out the
analysis by translation into an identity cell (see
Remark~\ref{r.unitcell}). However, it is more natural in this
context to consider directly general cells (not necessarily
containing the identity element). We shall take the latter route.

\begin{df}\label{d.cg.general} ($\C$-genericity for general cells)\\
Let $\C$ be a convex subset of $W$ containing the identity
element. A vector $f$ in the root space $V$ is {\rm $\C$-generic}
if:
\begin{itemize}
\item[{\bf (i)}] For all $t\in T_{\C}$, $\langle f,\al_t\rangle
\not\in \{0, 1, -1\}.$ \item[{\bf (ii)}] For all $t\in T_{\partial
\C}$, $\langle f,\al_t\rangle = \pm 1.$ \item[{\bf (iii)}] If
$w\in \C$, $s,t\in S$, $m(s,t)=3$ and $ws, wt\not\in \C$ then
$$
\varepsilon_{w,s}\langle f,\al_{wsw^{-1}}\rangle =
\varepsilon_{w,t} \langle f,\al_{wtw^{-1}}\rangle\;(= \pm1),
$$
where $\varepsilon_{w,s}:=1$ if and only if $\ell(ws)>\ell(w)$ and
$-1$ otherwise.
\end{itemize}
\end{df}

\begin{lem}\label{t.l3q11}
Let $q=1$, let $\C$ be a convex subset of $W$ containing the
identity element, let $\T:S \to E_{\C}$ satisfy Axiom $(B)$ and
the involution relation $ \T_s^2 = 1\qquad(\forall s\in S), $ and
assume that $ \ua_t, \ub_t, \db_t \ne 0\qquad(\forall t\in
T_{\C}). $ If $s,t\in T$ satisfy $m(s,t)=3$, then $\T_s\T_t\T_s =
\T_t\T_s\T_t$ if and only if:
\begin{itemize}
\item[(a)]
For every $w\in W$ such that $\ell(w)<\ell(ws)$, $\ell(w)<\ell(wt)$,\\
and $|w\langle s,t\rangle\cap \C|\ge 2$: $ {1\over
\ua_{wstsw^{-1}}} = {1\over \ua_{wsw^{-1}}} + {1\over
\ua_{wtw^{-1}}}. $ \item[(b)] For every $w\in \C$ such that
$|w\langle s,t\rangle\cap \C|= 1$,\ $ \oa_{wsw^{-1}} =
\oa_{wtw^{-1}} = \pm 1. $
\end{itemize}
\end{lem}

Proof is similar to the proof of Lemma~\ref{t.l3q1} and is
omitted.

\noindent{\bf Note:}
It follows that Theorems~\ref{t.f2rep}, \ref{t.rep2f},
\ref{t.f2rep.q} and \ref{t.rep2f.q} hold for general cells (with
Definition~\ref{d.cg} replaced by Definition~\ref{d.cg.general}).

\medskip

Recall the notation $\C_A(w)$ from Definition~\ref{d.A_equiv}.

\begin{df}\label{d.CF}
Let $w\in W$, and let $f$ be an arbitrary vector in the root space
$V$ of $W$. Let
$$
A = A_f := \{t\in T\,|\,\langle f,\al_t\rangle \in\{1, -1\}\}.
$$
\begin{itemize}
\item[(1)] Define
$$
\C^f(w):=\C_A(w).
$$
\item[(2)] If $f$ is $\C^f(w)$-generic then the corresponding AY
representation of $W$, with the $(SON)$ normalization (see
Subsection~\ref{s.normal}), will be denoted $\T^f(w)$.
\end{itemize}
\end{df}







For a finite simply laced Coxeter group $W$ let
$$
\delta := \frac{1}{2}\sum\limits_{t \in T} \alpha_t,
$$
half the sum of all the positive roots. ($\delta$ is more commonly
denoted $\rho$, but this letter has a different meaning in this
paper.) In this section we consider $f = \delta$. Recall also that
a finite crystallographic (in particular, simply laced) Coxeter
group is a Weyl group.

\begin{fac}\label{t.delta1}
For every reflection $t\in T$
$$
\langle \delta, \al_t\rangle=\hgt(\al_t),
$$
where the height $\hgt(\al_t)$ is as defined in
Section~\ref{s.prelim}.
\end{fac}

\noindent{\bf Proof.} $\delta$ is also equal to the sum of all
simple coroots (cf. \cite[\S 13.3, Lemma A]{Hum2}).

\qed

A characterization of the corresponding cells follows.

\begin{pro}\label{t.delta2}
For every finite simply laced Coxeter group $W$ and every element
$w\in W$, the cell $\C^\delta(w)$ is a standard (left) descent
class, i.e.,
$$
\C^\delta(w)=\{v\in W|\ \D_S(v)=\D_S(w)\}.
$$
\end{pro}

\noindent{\bf Proof.} By Fact~\ref{t.delta1}, $\langle\delta,
\al_t\rangle$ is a positive integer for every reflection $t\in T$,
and $\langle\delta, \al_t\rangle = 1$ if and only if $t\in S$.

\qed

\medskip

Note that $\delta$ is $\C^\delta(w)$-generic for every $w\in W$.

\begin{df}\label{d.des_rep}
Let $w\in W$. The representation $\T^\delta(w)$ of $W$ is called
the {\em descent representation} corresponding to $w$.
\end{df}

An immediate consequence is an analogue of the classical Young
Orthogonal Form
for this family of
representations.

\begin{cor}\label{t.yof1}{\bf (Orthogonal Form for Descent Representations)}\\
Let $W$ be a finite simply laced Coxeter group and let $w\in W$.
Then, for every Coxeter generator $s\in S$ and every element $v\in
\C^\delta(w)$, $\T:=\T^\delta(w)$ satisfies
$$
\T_s(C_v) = \cases{%
{1\over\hgt(\al_{vsv^{-1}})}C_v+\sqrt{1-{1\over\hgt(\al_{vsv^{-1}})^2}}
C_{vs}, & if $\ell(v)<\ell(vs)$;\cr
-{1\over\hgt(\al_{vsv^{-1}})}C_v+\sqrt{1-{1\over\hgt(\al_{vsv^{-1}})^2}}C_{vs},
& if $\ell(v)>\ell(vs)$.\cr}
$$
\end{cor}

This extends the classical Young Orthogonal Form from irreducible
$S_n$ representations 
to all descent representations of all finite simply laced Coxeter
groups.

A $q$-analogue follows easily from Theorem~\ref{t.f2rep.q}. An
extension to non simply-laced Weyl groups will be carried out
elsewhere.

\bigskip

Recall that a standard descent class is a union of Kazhdan-Lusztig
cells \cite[\S 7.15]{Hum}.

\begin{con}\label{t.descent}
For the $q$-analogue,
$$
\T^\delta(w) \cong \hbox{Kazhdan-Lusztig representation on\ }
\C^\delta(w)
$$
as ${\cal H}_q(W)$-representations.
\end{con}

Conjecture~\ref{t.descent} is related to the following general
problem.

\begin{prb}
Let $w_1, w_2\in W$ and let $f_1,f_2\in V$. When is
$$\T^{f_1}(w_1)\cong \T^{f_2}(w_2)\ \ \ ? $$
\end{prb}

This problem will be studied in \cite{U-II}.

\subsection{Irreducible Representations}

Let $\lambda$ be a partition of $n$, and let $Q$ be a standard Young
tableau of shape $\lambda$. For a permutation $\pi$ in the
symmetric group $S_n$ denote by $Q^\pi$ the tableau obtained
from $Q$ by replacing each entry $i$ by $\pi(i)$. Then

\begin{thm} The subset $\C_Q$ of the symmetric group $S_n$ given
by
$$
\C_Q:=\{\pi\in S_n|\ Q^{\pi^{-1}} \hbox{ is standard}\}
$$
is a minimal AY cell. Moreover, $\C_Q$ carries an AY representation
which is isomorphic to the irreducible Specht module $S^\lambda$.

In particular, every irreducible representation of the symmetric
group $S_n$ may be realized as a minimal abstract Young
representation.
\end{thm}

\noindent{\bf Proof (Sketch).} For $1\le i\le n$ let $c(i)$ be the
contents of $i$ in $Q$. Let $f_Q=(c(2)-c(1),\dots,c(n)-c(n-1))$ be
the hook distances vector of $Q$. One may verify that $f_Q$ is a
$\C_Q$ generic vector and that $\C^{f_Q}(id)=\C_Q$. Moreover, the
resulting representation matrices of the Coxeter generators are
given by the classical Young Form for the Specht modules.

\qed

A detailed proof, as well as a full characterization of AY pairs
in the symmetric group are given in~\cite{U-II}.

\bigskip

A similar theorem holds for Weyl groups of type $B$. Namely, every
irreducible representation of the hyperoctahedral group $B_n$ may
be realized as a minimal abstract Young representation. A detailed
proof will be given in~\cite{U-II}.

\section{Appendix 1: Equivalence of Axioms%
}\label{s.ax1}

The aim of this section is to prove Theorem~\ref{t.a1}, stating
the equivalence of Axioms $(A)$ and $(B)$ under suitable assumptions.

Recall that $(W,S)$ is a Coxeter system,
$(\T,\C)$ is a minimal AY pair for (the Iwahori-Hecke algebra of) $W$,
and
$$
a_s(w)=a_{s'}(w') \then b_s(w)=b_{s'}(w')
\qquad(\forall s,s'\in S, w,w'\in \C).
$$
%
%


In order to prove that $\T$ satisfies Axiom $(B)$
we first prove it for Coxeter groups with two generators, namely:
dihedral groups.

\medskip

\begin{lem}\label{t.conj.dihedral}
Let $(W,S)$ be a Coxeter system with two generators: $S = \{s_1,s_2\}$.
Let $(\T,\C)$ be a minimal AY pair for (the Hecke algebra of) $W$.
If $(w,s),(\tw,\ts)\in \C\times S$ satisfy
$\tw\ts\tw^{-1} = w s w^{-1}$
then
$$
a_{\ts}(\tw) = \cases{%
a_s(w), &if $\ell(ws)-\ell(w)$ and $\ell(\tw\ts)-\ell(\tw)$ have the same sign;\cr
a_s(ws), &otherwise.\cr}
$$
The first case occurs, for $(\tw,\ts) \ne (w,s)$, if and only if
$m = m(s_1,s_2) < \infty$ and
$(\tw,\ts) = (wsw_0,w_0sw_0)$, where $w_0$ is the longest element in $W$.
Note that
$$
w_0s_1w_0 = \cases{%
s_1, &if $m$ is even;\cr
s_2, &if $m$ is odd.}
$$
\end{lem}
\noindent{\bf Proof.}
First note that if $(\tw,\ts) = (ws,s)$ then clearly $\tw\ts\tw^{-1} = wsw^{-1}$,
but of course $\ell(ws)-\ell(w)$ and $\ell(\tw\ts)-\ell(\tw)$ have opposite signs.

If $m=\infty$ then $W = \langle s_1,s_2 \,|\, s_1^2 = s_2^2 =1\rangle$
is an infinite group, the free product of two groups of order $2$, and
its Cayley graph is a doubly infinite path. Every element of $W$ has
a unique reduced expression in terms of $s_1$ and $s_2$.
It follows that in this case
$$
\tw\ts\tw^{-1} = w s w^{-1} \iff
(\tw,\ts)\in \{(w,s), (ws,s)\},
$$
and our claim trivially holds.

Assume now that $2\le m<\infty$. $W$ is then the dihedral group of order $2m$,
and its Cayley graph is a $2m$-cycle.
Let $w_0$ be the unique longest element of $W$. It has two reduced expressions:
$$
w_0 = s_1 s_2 \cdots = s_2 s_1 \cdots,
$$
where both products have length $m$ and their factors alternate between
$s_1$ and $s_2$. Every other element of $W$ has a unique reduced expression,
of length at most $m-1$.
$W$ contains exactly $m$ reflections, and each of them has the form $wsw^{-1}$
for some $w\in W$ and $s\in S$. It can also be written as $(ws)s(ws)^{-1}$,
$(ww_0)(w_0sw_0)(ww_0)^{-1}$ and $(wsw_0)(w_0sw_0)(wsw_0)^{-1}$.
These $4$ representations for a reflection are distinct,
since the elements $w$, $ws$, $ww_0$ and $wsw_0$ are distinct.
We thus obtain at least $4 \cdot m$ representations for reflections.
The set of pairs $W \times S$ has $2m \cdot 2 = 4m$ elements, and therefore
each reflection is obtained {\em exactly} $4$ times in the form $wsw^{-1}$.
Also, $\ell(ww_0) = m-\ell(w)$ so that
$$
\ell(w)<\ell(ws) \iff \ell(wsw_0) < \ell(ww_0)=\ell(wsw_0 \cdot w_0sw_0).
$$
Returning to our setting, we only need to show that if $w, wsw_0\in\C$ then
$$
a_s(w) = a_{w_0 s w_0}(wsw_0).
$$

Assume for concreteness that $s = s_1$.
If $w, ws_1w_0\in\C$ then there is a unique shortest path
(of length $m-1$) from $w$ to $ws_1w_0$ in the Cayley graph of $W$:
$$
w \to ws_2 \to ws_2s_1 \to \cdots \to ws_1w_0.
$$
$\C$ is a minimal AY cell, and is therefore convex
(by Corollary~\ref{t.MAY_cell}(i)), so all the vertices in this path
belong to $\C$. By Corollary~\ref{t.MAY_cell}(ii), all the edges of the path
carry nonzero coefficients $b_.(.)$. This path corresponds to the unique
reduced expression (of length $m-1$) for $s_1w_0$.
This expression is obtained uniquely from each of the two reduced expressions
for $w_0$, by deleting one letter: the first letter from $s_1s_2\cdots$,
and the last letter from $s_2s_1\cdots$.
Comparing the coefficients of $C_{ws_1w_0}$ in the two expressions for
$\T_{w_0}(C_w)$,
$$
\cdots \T_{s_2}\T_{s_1}(C_w) = \cdots \T_{s_1}\T_{s_2}(C_w),
$$
shows that
$$
a_{s_1}(w) \cdot \Pi(b) = \Pi(b) \cdot a_{w_0 s_1 w_0}(w s_1 w_0)
$$
where
$$
\Pi(b) := b_{s_2}(w) b_{s_1}(ws_2) \cdots \ne 0.
$$
Therefore
$$
a_{s_1}(w) = a_{w_0 s_1 w_0}(w s_1 w_0),
$$
as claimed.

\qed

For the general case of Theorem~\ref{t.a1} we need a connectivity argument,
supplied by the following result.

\begin{lem}\label{t.conj.seq}
Let $(W,S)$ be a Coxeter system, and let $(w,s),(\tw,\ts)\in W \times S$.
Then
$$
\tw\ts\tw^{-1}= w s w^{-1}
$$
if and only if there exists a sequence of pairs
$$
(w,s)\stackrel{\varepsilon}{\too}
(w_1,s_1)\stackrel{\bs_1}{\too}
(w_2,s_2)\stackrel{\bs_2}{\too}
\ldots\stackrel{\bs_{k-1}}{\too}
(w_k,s_k) =
(\tw,\ts)
$$
such that:
\begin{enumerate}
\item
For $1\le i\le k$:
$(w_i,s_i)\in W \times S$,
$w_i s_i w_i^{-1} = w s w^{-1}$, and
$\tw^{-1}w_i < \tw^{-1}w_i s_i$.
\item
$s_1 = s$ and $w_1 = ws^{\varepsilon}$, where
$$
\varepsilon := \cases{%
0,&if $\tw^{-1}w < \tw^{-1}ws$;\cr
1,&if $\tw^{-1}w > \tw^{-1}ws$.\cr}
$$
\item
For $1\le i\le k-1$:
$\bs_i\in S \setminus \{s_i\}$, and if $(s_i \bs_i)^{m_i} = 1$ is
the corresponding Coxeter relation $(m_i\ge 2)$ then $m_i\ne\infty$ and
$$
s_{i+1} = \cases{%
s_i, &if $m_i$ is even;\cr
\bs_i, &if $m_i$ is odd.\cr}
$$
Also,
$$
w_i^{-1}w_{i+1} = \bs_i s_i \cdots,
$$
a product of length $m_i-1$, with factors alternating between $\bs_i$ and $s_i$.
\item
$$
\tw^{-1}w \ge \tw^{-1}w_1 > \tw^{-1}w_2 > \ldots > \tw^{-1}w_k =1
$$
(in right weak Bruhat order), and
$$
\ell(\tw^{-1}w_i) - \ell(\tw^{-1}w_{i+1}) = m_i - 1 \qquad(1\le i\le k-1).
$$
\end{enumerate}
\end{lem}

\noindent
The lemma gives an explicit reduced expression:
$$
w^{-1}\tw = s^{\varepsilon} (\bs_1 s_1 \cdots) \cdots (\bs_{k-1} s_{k-1} \cdots).
$$
Note also that $s_{i+1} \ne s_i$ if and only if $m_i$ is odd, and then
$m_i = m(s_i,s_{i+1})$.

\begin{cor}\label{t.conj}{\rm \cite[Ch. IV, \S 1, Prop. 3]{B}}
Two simple reflections $s, \ts \in S$ are conjugate in $W$ if and only if they are
connected, in the Dynkin diagram of $W$, by a path consisting entirely of edges with
odd labels.
\end{cor}

\noindent{\bf Proof of Lemma~\ref{t.conj.seq}.}
Clearly, if such a sequence exists then $\tw\ts\tw^{-1} = w s w^{-1}$.

Conversely, assume that $\tw\ts\tw^{-1} = w s w^{-1}$.
Let $s_1:=s$ and $w_1:=ws^{\varepsilon}$, where
$$
\varepsilon := \cases{%
0,&if $\tw^{-1}w < \tw^{-1}ws$;\cr
1,&if $\tw^{-1}w > \tw^{-1}ws$.\cr}
$$
Then clearly $w_1 s_1 w_1^{-1} = wsw^{-1}$, $\tw^{-1}w_1 < \tw^{-1}w_1 s_1$
and $\tw^{-1}w \ge \tw^{-1}w_1$.

Inductively, assume that $\tw^{-1}w_i < \tw^{-1}w_i s_i = \ts\tw^{-1}w_i$.
If $\tw^{-1}w_i = 1$ then $(w_i,s_i) = (\tw,\ts)$ and we are done.
Thus assume that $\tw^{-1}w_i \ne 1$. Choosing a reduced expression for
$\tw^{-1}w_i$ gives two distinct reduced expressions for
$\tw^{-1}w_i s_i = \ts\tw^{-1}w_i$.
One of these has $s_i$ as the rightmost letter, while the other does not.
Any reduced expression can be transformed into any other reduced expression
of the same group element by a sequence of braid moves
\cite[Theorem 3.3.1(ii)]{BB}.
In our case, at some point during this process $s_i$ ceases to be
the rightmost letter in the
expression. This means that at this point we use a braid relation
$\cdots \bs_i s_i = \cdots s_i \bs_i$ ($m_i\ge 2$ letters on each side),
where $\bs_i \in S\setminus \{s_i\}$ and $(s_i \bs_i)^{m_i} = 1$.
Thus there exists a reduced expression $r_i$ for which
\begin{equation}\label{sofsof}
\tw^{-1}w_i s_i = r_i \cdots \bs_i s_i = r_i \cdots s_i \bs_i
\end{equation}
(two reduced expressions), where the length of each ``$\cdots$" is $m_i-2\ge 0$.
In particular, $\tw^{-1}w_i = r_i \cdots \bs_i$ so that
$$
\ts r_i \cdots \bs_i = \ts\tw^{-1}w_i = \tw^{-1}w_i s_i =
r_i \cdots \bs_i s_i = r_i \cdots s_i \bs_i.
$$
The ``tails'' (of length $m_i-1$) in the first and last expressions
in the above line are equal.
Thus either $\ts r_i = r_i s_i$ (for $m_i$ even)
or $\ts r_i = r_i \bs_i$ (for $m_i$ odd).
Let $w_{i+1} := \tw r_i$ and
$$
s_{i+1} := \cases{%
s_i &if $m_i$ is even;\cr
\bs_i &if $m_i$ is odd.\cr}
$$
Then $w_{i+1}s_{i+1} = \tw r_i s_{i+1} = \tw\ts r_i$, so that
$w_{i+1}s_{i+1}w_{i+1}^{-1} = \tw\ts\tw^{-1}$.
Also $\tw^{-1}w_{i+1} = r_i < r_i s_{i+1} = \tw^{-1}w_{i+1}s_{i+1}$
and $w_i^{-1}w_{i+1} = \bs_i \cdots$ (of length $m_i-1$).

Since $\ell(\tw^{-1}w_i) - \ell(\tw^{-1}w_{i+1}) = m_i - 1 > 0$
by (\ref{sofsof}), this process must stop
after a finite number of steps, and then necessarily $\tw^{-1}w_k = 1$ as above.
This completes the proof.

\qed

\bigskip

\noindent
{\bf Proof of Theorem~\ref{t.a1}.}
We assume that $\T$ satisfies Axiom $(A)$ and that
$(\T,\C)$ is a minimal AY pair.
As noted above, we only need to show that if $(w,s),(\tw,\ts)\in \C\times S$
satisfy $\tw\ts\tw^{-1} = wsw^{-1}$ and, say, $\ell(w) < \ell(ws)$ as well as
$\ell(\tw) < \ell(\tw\ts)$, then $a_s(w) = a_{\ts}(\tw)$.

By Lemma~\ref{t.conj.seq}
there exists a geodesic (with special structure) from $w$ to $\tw$ in $W$.
By the convexity of $\C$ (Corollary~\ref{t.MAY_cell}(i)),
this geodesic is actually in $\C$.
This geodesic decomposes into segments, each of which belongs to a coset
in $W$ of a parabolic subgroup with two generators, $s_i$ and $\bs_i$.
Lemma~\ref{t.conj.dihedral} and the special properties of the path
complete the proof.

\qed


\section{Appendix 2: The $b$-independence of Characters}\label{s.ax2}

The aim of this section is to prove the following result,
answering affirmatively Problem~\ref{p.char}:

\begin{thm}\label{t.b_indep}
Let $(\T,\C)$ be an AY pair satisfying Axiom $(B)$. If $\C$ is finite then
all the corresponding character values are polynomials in $a$-coefficients
only (no $b$-coefficients).
\end{thm}

We shall need the following result.

\begin{lem}\label{t.disconnect}
Fix a reflection $t\in T$. Deleting from the Cayley graph
$\G(W,S)$ all the edges $\{w,ws\}$ for which $wsw^{-1} = t$ leaves
exactly two connected components. Each deleted edge has one vertex
in each of these components.
\end{lem}
\noindent{\bf Proof.} By Corollary~\ref{t.equiv_gdc} (with $A =
\{t\}$), the connected components are exactly the nonempty
generalized descent classes of the form $W_{\{t\}}^D$, where
either $D = \emptyset$ or $D = \{t\}$. These two classes are
clearly nonempty and, by Theorem~\ref{t.tits}, convex; thus, in
particular, connected. The remaining claim follows from
Proposition~\ref{t.des_diff}.

\qed

\noindent{\bf Proof of Theorem~\ref{t.b_indep}.} Fix $g\in W$ and
write it as a product of generators:
$$
g = s_{i_k} \cdots s_{i_1}\qquad(s_{i_j}\in S,\,\forall j).
$$
Let $\langle \cdot, \cdot \rangle$ be the inner product on
$V_{\C}$ for which $\{C_w \,|\, w\in \C\}$ is an orthonormal
basis. The character value at $g$,
$$
\tr(\T_g) = \sum_{w\in \C} \langle \T_g(C_w),C_w \rangle
          = \sum_{w\in \C}
            \langle \T_{s_{i_k}} \cdots \T_{s_{i_1}}(C_w),C_w \rangle,
$$
is, according to Axiom $(A)$, a sum of products of the form
\begin{equation}\label{e.prod}
x_1(w_0) x_2(w_1) \cdots x_k(w_{k-1}),
\end{equation}
where each $x_j$ is either $a_{s_{i_j}}$ or $b_{s_{i_j}}$ ($1\le
j\le k$),
\begin{eqnarray*}
w_0 &:=& w,\\
w_j &:=& \cases{w_{j-1}, &if $x_j = a_{s_{i_j}}$\cr
                w_{j-1}s_{i_j}, &if $x_j = b_{s_{i_j}}$\cr}
         \qquad(1\le j\le k),
\end{eqnarray*}
and the sum is over all the choices of $x_1,\ldots,x_k$ for which
$w_k=w_0$ (and $w_j\in\C$, $\forall j$). It suffices to show that
each of these products can be rewritten in terms of
$a$-coefficients only.

Consider a fixed product~(\ref{e.prod}), with $w_k = w_0$.
Geometrically, this product corresponds to a closed walk in the
Cayley graph of $W$: $w_0 \to w_1 \to \ldots \to w_k = w_0$. Each
step in this walk is either a loop $w_{j-1} \to w_{j-1}$ (if $x_j
= a_{s_{i_j}}$) or a directed edge $w_{j-1} \to w_{j-1}s_{i_j}$
(if $x_j = b_{s_{i_j}}$). Now, according to Axiom $(B)$,
$$
b_s(w) \in \{\ub_{wsw^{-1}}, \db_{wsw^{-1}}\}.
$$
It suffices to show that all the $b$-coefficients in the
product~(\ref{e.prod}) can be paired into products of the form
$\ub_t \db_t$, which by Lemma~\ref{t.l1}(a) can be expressed using
$a$-coefficients only.

Indeed, fix $t\in T$. By Lemma~\ref{t.disconnect}, the edges $\{w,ws\}$
for which $wsw^{-1}=t$ occur in any given closed walk an even number of
times. Moreover, exactly half of these occurrences correspond to
factors $\ub_t$, and half to $\db_t$ (by Axiom $(B)$,
Proposition~\ref{t.des_diff}, and the characterization of the
connected components in the proof of Lemma~\ref{t.disconnect}).
This completes the proof.

\qed

\end{document}